\documentclass[nonblindrev]{informs3}

\OneAndAHalfSpacedXI
\usepackage{graphicx}
\usepackage[utf8]{inputenc}
\usepackage[table]{xcolor}
\usepackage{amsmath,amssymb} %amsthm
\usepackage{multirow}
\usepackage{authblk}
\usepackage{url}
\usepackage{pgf}
\usepackage{color}
\usepackage{todonotes}
\usepackage{booktabs} %for table weighted rules
\usepackage{algorithm}
\usepackage{algorithmic} 
\usepackage{pgfplots}
\usepackage{soul}
\usepackage{bbm}
\usepackage{hyperref}
\usepackage{enumerate}

\usepackage{tikz-network}

\newtheorem{prop}{Proposition}

\newtheorem{remark}{Remark}

\newtheorem{example}{Example}

\pgfplotsset{compat=1.17}
\usepackage{natbib}
 \bibpunct[, ]{(}{)}{,}{a}{}{,}%
 \def\bibfont{\small}%

\usepackage{cleveref}
\usepackage{geometry}
\usepackage{setspace}
\allowdisplaybreaks
\newcommand{\esra}{\textcolor{blue}}
\newcommand{\ignore}[1]{}

\TITLE{Optimizing Strategic and Operational Decisions of Car Sharing Systems under Demand Uncertainty and Substitution} %Designing Car Sharing Systems with Electric Vehicles under Demand/Adoption Uncertainties with Quality of Service Constraints
\ARTICLEAUTHORS{
\AUTHOR{Sinan Emre Kosunda} \AFF{Industrial Engineering Program, Sabanc{\i} University, Istanbul, Turkey, \EMAIL{sinankosunda@sabanciuniv.edu}}
  \AUTHOR{Beste Basciftci\thanks{Corresponding author}} \AFF{Department of Business Analytics, Tippie College of Business, University of Iowa, Iowa City, Iowa, USA, \EMAIL{beste-basciftci@uiowa.edu}}
  \AUTHOR{Esra Koca} \AFF{Industrial Engineering Program, Sabanc{\i} University, Istanbul, Turkey, \EMAIL{esra.koca@sabanciuniv.edu}}
}

\begin{document}

\ABSTRACT{
%\textbf{\textit{Problem Definition:}} 
Optimizing car sharing systems under demand uncertainty is an emerging problem for ensuring profitable and sustainable operations of these services while taking into account quality of service concerns. With the increasing adoption of electric vehicles and environmental awareness, this problem requires consideration of a mix fleet of vehicles with gasoline-powered and electric, complicating the strategic and operational planning as the demand of each vehicle type are observed. 
%\textbf{\textit{Methodology/Results:}}
To address this problem, we propose a two-stage stochastic mixed-integer program leveraging spatial-temporal networks that  capture the strategic and operational decisions of these systems over a multi-period planning horizon. We optimize the location decisions of regions to serve with purchasing decisions of the vehicles in the first-stage problem under budget and carbon emission considerations in designing the fleet, while considering parking capacities, satisfying one-way and round-trip car rental requests, and relocating cars between open regions under each demand realization in the second-stage problem. 
We then introduce demand substitution to this problem by extending and generalizing the multi-commodity formulation, and allowing satisfaction of customer demand of each vehicle type with its alternatives. We further prove that the corresponding second-stage problem has a totally unimodular constraint matrix. By benefiting from this result, as our solution approach, we provide a branch-and-cut based decomposition algorithm with enhancements. 
We present an extensive computational study highlighting the value of the proposed models from different perspectives and demonstrating the performance of the proposed solution algorithm with significant speedups.
%\textbf{\textit{Managerial implications:}}
Our case study demonstrates the benefits of incorporating strategic and operational decisions with the integration of demand substitution, and provides insights for region opening and fleet allocation plans under demand uncertainty. 
Furthermore, introducing substitution to the car sharing operations leads to higher quality of service and flexibility in operations with lower costs under various settings. 
}


\maketitle

\section{Introduction}

Car sharing has been an emerging area of smart city operations by utilizing vehicles in a more efficient manner while reducing congestion and providing environmental benefits with the help of its sharing ecosystem \citep{CarSharingBenefits2016}. 
In car sharing services, customers rent cars for a certain amount of time, where they pick up and drop off the vehicles at regions that are served by these service provider companies. Enabling pick up and drop off locations to be different from each other provides flexibility to the customers by allowing one-way trips in addition to the round-trips, in which these locations are the same. Nevertheless, this necessitates more complex operational planning, as the car sharing companies need to ensure rebalancing of vehicles between different service regions throughout their daily operational plans. Another complexity arises when the customer demand is not fully known in advance. Furthermore, different customer segments can prefer different vehicle types, leading to consideration of various vehicle types simultaneously in fleet management and allocation.  
Consequently, the car sharing companies need to determine their service regions while taking into account a mix fleet of vehicles and their operational planning under demand uncertainty. 

Transportation sector constitutes the largest source of emissions of carbon dioxide in the United States \citep{transportationEmission}. 
Thus, increasing usage of electric vehicles (EVs) within the mobility systems are projected to significantly alleviate the impact of emissions over the next decade. 
Despite of their fuel-efficiency and environmental benefits, the major concerns against the mass adoption of EVs are the higher purchasing costs compared to gasoline vehicles, range anxiety and limited charging infrastructure. 
Therefore, instead of the ownership of EVs, car sharing companies can provide rental services to the customers through these vehicles. 
This approach can further eliminate the concerns of the customers by removing the higher cost of purchasing and providing charging infrastructure within the parking stations operated by the car sharing companies \citep{Brandstatter2020}. 
As car sharing services are mainly operating in urban areas, range anxiety can be mitigated with easier access to the charging infrastructure. 
Nevertheless, with the improvements in the battery technologies, most of the recently sold EVs in the United States in 2022 have longer range with more than 240 miles \citep{EVSales2022}. Such a driving range is significantly larger than the distance travelled in majority of the car sharing trips, considering the car sharing operations of various companies in the United States \citep{He2017} and Germany \citep{strohle2019}. 
These factors can eliminate the need of recharging of EVs throughout the day, which can make them more attractive to the customers and easier to operate for the car sharing companies. 

%bu paragrafa kadar hic we do that ... tarzi birsey yazmamistik. Burada carbon emission kisitindan oyle bahsettim. Bir sekilde emisyon kisitindan introda da bahsetmek lazim ama daha farkli bahsedilebilir mi?
Although EVs are becoming an important component of the car sharing services, their operations need to be jointly considered with the gasoline  vehicles as these companies can have existing fleets with gasoline vehicles and adoption of EVs by customers can require a transitional period \citep{AboueeMehrizi2021}. 
For instance, the largest car sharing companies, such as Zipcar and Share Now, operate a mix fleet with electric and gasoline vehicles. Thus, the joint consideration of both vehicle types constitute a critical concern during this transition. Moreover, since the purchasing cost of EVs is higher than the gasoline vehicles, environmental consciousness of the companies need to be explicitly captured while constructing the mixed fleet of electric and gasoline vehicles. %there should be a motivation for the sharing companies to buy more EVs to incorporate their environmental consciousness.
\textcolor{black}{Therefore, to increase the percentage of the EVs in the fleet, it becomes critical to restrict carbon emissions to guarantee that the average carbon emission of the fleet is limited by a given unit carbon emission allowance.} %which is in line with the most commonly used carbon emission calculation protocol \citep{greenhousegas}, that is considered in production planning problems \citep{absi2013, absi2016}, but has not been utilized in the literature in optimizing the car sharing systems.}
\iffalse
Therefore, to increase the percentage of the EVs in the fleet, we consider carbon emission constraints that should be taken into account,  %. Instead of considering the actual carbon emissions of the trips \citep{Chang2017}, we 
which guarantee that the average carbon emission of the fleet is limited by a given unit carbon emission allowance. Although a similar approach, which is adopted from the most commonly used carbon emission calculation protocol \citep{greenhousegas}, is considered in production planning problems \citep{absi2013, absi2016}, it has not been utilized in the literature in optimizing the car sharing systems. 
\fi

As car sharing companies need to consider a mix fleet of vehicles, the demand for each vehicle type need to be incorporated into the planning of service region design, operations and repositioning activities. To this end, leveraging substitution between different vehicle types in satisfying customer demand provides flexibility in operations with higher quality of service and results in higher revenues. 
In addition to the vehicles with different fuel technologies, mix fleets can be classified in terms of the vehicle segments, and substitution between different vehicle types can become further favorable in satisfying the customer demand corresponding to each segment. 
Although the value of substitution has been demonstrated across various operations management problems in efficiently matching supply and demand, its value and integration to car sharing operations have not been investigated. 
To address these problems, we study the strategic and operational decision making of the car sharing companies with a mix fleet of vehicles by capturing the uncertainty in customer demand and integration of substitution between different vehicle types. The goal of the companies is to maximize their net profits over this multi-period planning problem, while taking into account the revenue obtained from round-trip and one-way rentals and the cost of serving regions and rebalancing vehicles between them for satisfying customer demand. 
Consequently, our contributions can be summarized as follows:
\begin{itemize}
    \item We propose a service region design and operational planning problem for a car sharing company with a mix fleet of vehicles by using spatial-temporal networks to model the round-trip and one-way rentals along with the repositioning trips. To integrate uncertainty in customer demand, we formulate this problem as a two-stage stochastic mixed-integer program that maximizes the net profit of the company, where the first-stage problem determines which regions to serve and how to allocate the fleet of each vehicle type to each region under a budget limitation and carbon emission considerations. Given these decisions, the second-stage problem optimizes the operational plans under each demand realization. 
    \item We introduce substitution to this problem by allowing customer demand of each vehicle type to be satisfied by its potential alternatives, with a penalty incurred by the car sharing company from its regular prices for incentivizing such trips. To formulate this problem, we propose a two-stage stochastic mixed-integer program that extends our spatial-temporal network representation to capture the substitute trips for round-trip and one-way rentals for each vehicle type, and show that it generalizes our initial formulation. 
    \item We develop an exact decomposition based algorithm to solve the resulting challenging problems, where we leverage the structure of the proposed formulations. We further provide computational enhancements to improve the efficiency of the algorithm. Our computational study illustrates significant speed-ups obtained by our solution algorithm with enhancements in comparison to the off-the-shelf solvers. 
    \item We present an extensive case study based on real data sets. Our findings demonstrate that introducing substitution to the car sharing operations significantly impacts service region design, fleet management decisions, demand satisfaction and operational plans. We provide various managerial insights that can be summarized as follows:
    \begin{itemize}
    \item When substitution is allowed between different vehicle types, customer demand can be satisfied more, leading to higher quality of service and higher profits. Depending on the demand distribution across different service regions, this can further impact which regions to serve. 
        \item When the budget of the company for constructing its fleet is smaller, the impact of substitution becomes larger, leading to higher percentage improvements in net profit in comparison to the plans without substitution. 
        \item With substitution, car sharing companies need to perform less relocation trips between different service regions over time, since one-way trips can be utilized more to satisfy the necessary rebalancing operations. 
        \item Depending on the penalty amount of substitution incurred by the car sharing company, amount of rental trips that are satisfied by alternative vehicle types changes, where the smaller penalty values leads to higher substitution rates.
        \item Carbon emission limitations considered affect the fleet allocation to each vehicle type, such as electric and gasoline vehicles, impacting the value of substitution remarkably.  
    \end{itemize} 
\end{itemize}

The remainder of the paper is organized as follows. In Section \ref{sec:Literature}, we provide the relevant literature. In Section \ref{sec:ProblemForm}, we first present the problem statement with the spatial-temporal network, and then formulate the service region and operational planning problems by introducing substitution. 
In Section \ref{sec:solution}, we present our solution algorithm, and in Section \ref{sec:ComputationalStudy}, we provide our detailed case study with various insights, showcasing the importance of the presented models and the computational efficiency of the solution algorithm. Section \ref{sec:Conclusion} concludes the paper with final remarks. 

\section{Literature Review}
\label{sec:Literature}
The demand for car sharing services is increasing rapidly due to the recent economical and environmental circumstances, and the constant growth of transportation requirements in daily lives, which brings new problems for car sharing companies. Excessive competition and increasing costs also make companies to deal with the strategic and operational decisions more seriously. Thus, the literature on various aspects of car sharing systems is growing fast \citep{ferrero2018, nansubuga2021}.

The problems observed in car sharing systems are studied from different perspectives under different assumptions and modeling approaches. Earlier studies on car sharing systems consider a deterministic setting where all problem parameters are known beforehand \citep[e.g.][]{CORREIA2012,NAIR201447,BOYACI2015718,Chang2017,GAMBELLA2018234} while recent studies are mostly focused on stochastic settings where some problem parameters, i.e. the demand, are not known with certainty \citep[e.g.][]{kaspi2014,He2017,Luetal2017,CALIK2019, Basciftci2019_DDDR,Zhang2020}. The latter studies also differ from each other based on the modeling scheme used for dealing with uncertainty. Markovian models \citep[e.g.][]{kaspi2014}, robust optimization \citep[e.g.][]{He2017} and two-stage stochastic programming \citep[e.g.][]{Luetal2017} are the three main modeling approaches used in these studies. Besides, there are many studies using simulation to evaluate different strategical and/or operational strategies for car sharing systems \citep[e.g.][]{barth1999simulation, kaspi2014, pfrommer2014dynamic,jorge2014comparing}.

Since location decisions are expensive strategic decisions, they should be carefully made by the firms to optimize their performances and operations over a long term. We refer the interested reader to the book by \cite{locationscience} for a review of different location problems and their application areas. Locating the stations and fleet sizing and positioning are the two main decisions considered in the literature on car sharing systems. A variety of papers in the car sharing literature consider both of these decisions %(plus additional operational decisions including the accepted/rejected reservations, relocation, charging EVs, etc.) 
\citep[e.g.][]{CORREIA2012,BOYACI2015718,He2017,CALIK2019,Zhang2020}, however a fewer of these studies consider the additional operational decisions including the accepted/rejected reservations, relocation, charging EVs, etc. On the other hand, there are studies that focus on the operational decisions and optimize only the vehicle flows between the stations to maximize the profit \citep[e.g.][]{kaspi2014,pfrommer2014dynamic,Chang2017,Luetal2017,GAMBELLA2018234}. 

Other aspects that create a distinction between the studies in the car sharing literature are the types of the vehicles (identical or different) and the types of the trips allowed (one-way, round-trip or mixed) in the system. Although car sharing companies have several different car types in their fleets in practice, there are only a very few studies that consider a mix fleet of vehicles \citep{Chang2017, AboueeMehrizi2021}. To the best of our knowledge, \cite{AboueeMehrizi2021} is the only study that considers a mixed fleet with a focus on the interplay between electric and traditional vehicles. The authors model the problem 
as a queuing network and derive conditions under which it is optimal to use EVs in the system. As pointed out by \cite{Chang2017}, one-way trips give a great flexibility to the customers for planning their trips, but managing the system with one-way trips is a harder and costly problem for the companies due to the supply-demand imbalance across the stations caused by one-way trips. Accordingly, studies considering one-way trips mostly focus on the relocation actions for mitigating this imbalance - see the review by \cite{illgen2019literature} for the studies on relocation problems arising in one-way car sharing systems. \cite{ferrero2018} states that almost 50\% of the papers (among 137 papers published between 2001 and 2016) in the car sharing literature focus on only one-way trips. However, in practice, the companies allowing one-way trips also offer round trips with slightly less charges. For instance, ZipCar offers a new service ZipCar Flex in UK where both one-way and round-trips are allowed. Hence, considering both trip types in the system makes the study more realistic. 
%Note that allowing only round-trip rentals in the system makes the problem relatively easier to solve with a simpler operational problem. %and a detailed operational problem considering spatio-temporal rebalancing and substitution actions along with strategic level decisions are omitted. 

Increasing customers awareness and the successful performance increases in EVs, encourage car sharing companies to include EVs in their fleets. Accordingly, the research on green car sharing problems has grown recently \citep{ferrero2018}. When the fleet includes EVs, additional operations such as charging decisions for EVs and locating charging stations might be also considered in the problem \cite[e.g.][]{Zhang2020, He2021}. Besides, the affects of car sharing systems on carbon emissions are analyzed by different techniques \citep{Chang2017,jung2018emission,amatuni2020emission,luna2020emission}. However, to the best of our knowledge, the carbon emission constraints are considered in the design stage of a car sharing system only in \cite{Chang2017}, where %The carbon emission constraints of \cite{Chang2017} 
these constraints consider the actual flow of vehicles and calculate the exact emission of the system. However, the common approach for determining the emission for a system, product, etc. is to determine the average carbon emission per unit used or produced \citep{greenhousegas, absi2013}, and we adapt this approach in our study.

Substitution is widely used in manufacturing systems in different ways such as downward substitution or both way substitution \citep[e.g.][]{Rao2004, Shumsky2009, Dawande2010, Xu2011, LIU2019999, Feng2022}. It is a common approach used in practice since it works for the benefits of both sides. When substitution is allowed, customers might satisfy their requirements by a different product type offered instead of a shortage. Hence, substitution might increase customers' utility. Besides, substitution gives some flexibility to the producer to hedge against uncertainty of the demand when there exist different product types or raw materials in the system. \cite{Feng2022} in their study on online retailing point out substitution in car sharing as a possible interesting future research problem. To the best of our knowledge, \cite{SMET2021107703} is the only study that considers substitution in a car sharing context. The authors focus on %vehicle substitution in a car sharing 
a car sharing system with only round-trip rentals and develop a two-stage stochastic model to maximize the expected profit under uncertain demand. They allocate the fleet to the stations in the first stage and determine the requests that will be accepted and rejected along with substitution decisions in the second stage. The authors consider a setting where only downward substitution is allowed, i.e. higher segment cars can be used to satisfy the demand for a lower segment car.  
\textcolor{black}{Thus, authors consider a simpler setting by only allowing one-way substitution while disregarding service region and fleet design decisions, one-way rental trips and detailed consideration of operational problem with relocations.}

\textcolor{black}{In this study, we consider the service region design, fleet sizing and spatio-temporal operational planning problem of a car sharing system with a mix fleet of vehicles under uncertain demand with budget and carbon emission considerations while satisfying both round-trip and one-way rentals. We introduce substitution to this problem, complicating the operational problem in return for enhancing  the quality of service and reducing the operational costs. We develop a cutting-plane based solution algorithm tailored for this problem, and provide managerial insights on our case study demonstrating the value of substitution and optimization of strategic and operational decisions in car sharing systems.}

%bu kismi sana biraktim Beste, dokunmadim:
\ignore{
Our motivation is to consider an operational problem for a mixed fleet of vehicles under uncertain demand for both one-way and round-trip rentals with strategic region opening and fleet sizing decisions. To this end, we develop a mixed integer stochastic programming model. In the first stage of which region opening, fleet size of each car type with respect to the given budget and carbon emission rates, and the allocation of this fleet to service zones are decided. In the second stage we construct a spatial-temporal network for vehicle movements in each scenario to represent revenue of uncertain demand for one-way and round-trip rentals. Furthermore, we define substitution among each type of car to increase car sharing company’s revenue by enhancing demand satisfaction levels and decreasing demand loss. To increase computational efficiency, we develop a cutting-plane based solution algorithm tailored for this problem. We study the effect of different parameter settings for carbon emission and penalty rates, budget values, and value of substitution for managerial insights.   
}
\section{Problem Formulations}
\label{sec:ProblemForm}

In this section, we formally introduce the service region and operational planning problem for optimizing strategic and operational decisions of car sharing systems with a mix fleet of cars under demand uncertainty and car type substitution. We first present the problem statement with the spatio-temporal network that is used for capturing the operational level decisions in Section \ref{sec:ProblemStatement}. Then, we propose the service region and operational planning problem under demand uncertainty in Section \ref{sec:SROP}. Finally, we generalize this problem and introduce car type substitution in Section \ref{sec:SROP-S} by allowing demand of each car type to be satisfied by its alternatives. 

\subsection{Problem Statement and Spatial-Temporal Network}
\label{sec:ProblemStatement}

We consider a car sharing company that plans its service regions and fleet sizes for its mix fleet of cars while taking into account operational decisions under demand uncertainty. Customer demand can be identified through reservation of one-way trips and round-trips, where the one-way trips allow customers to pick up and drop off their cars at different service regions and the round-trips require customers to drop off their cars to their pick up location. The operational decisions are based on the car movements to satisfy one-way and round-trip customer demand and relocate cars when it is necessary. 
In particular, relocation trips are conducted by the car sharing company to rebalance the cars from one service region to another depending on the fleet allocation and customer demand.  

The goal of the car sharing company is to maximize its annual profit by considering the revenue obtained from the one-way and round-trips along with the cost of relocating cars and operating service regions. 
Company has a budget limitation in constructing its fleet, which consists of a mixture of car types to satisfy needs of different customer groups. To this end, the cars can be classified in terms of the customer market segments they will be preferred by or in terms of its fueling technology such as gasoline and electric vehicles. 
For this study, we focus on an environmentally-aware company that aims to design its fleet while taking into account its carbon emissions. Thus, purchase of gasoline cars needs to be adjusted with the purchase of electric cars to satisfy the emission targets for the fleet. 

We formulate this problem through two-stage stochastic mixed-integer programs, where the first-stage problems determine the strategic decisions including which regions to serve and how to construct the fleet and allocate to these regions. Given these decisions, the second-stage problems optimize the operational plans corresponding to the car movements resulting from one-way, round-trip and relocation trips under each demand realization. We approach this problem through two different formulations where the initial formulation does not allow substitution of demand through different car types and the latter formulation introduces substitution to increase profitability and customer satisfaction. Both formulations share the same set of strategic level decisions, resulting in the same first-stage problem. However, the second-stage problems of these formulations differ as allowing substitution complicates the problem significantly, requiring development of alternative formulations. To capture the operational decisions, both formulations benefit from spatial-temporal networks, which consider car movements on different service regions over a multi-period planning horizon. 

In the first-stage problem, the car sharing company determines which regions to serve from the set of possible service regions, denoted by $I$. The binary variable $z_i$ indicates whether region $i \in I$ is opened or not. To construct its fleet and allocate them to the open regions, the company considers car types from the set $K$. The integer variable $x_{ik}$ indicates the number of type $k \in K$ cars allocated to region $i \in I$ at the beginning of the planning. To open and operate a service region $i \in I$, the company needs to pay a fixed cost of $f_i$. Furthermore, company has a budget of $B$ for purchasing the cars, where each car type $k \in K$ has a cost of $c_k$ and emission amount of $e_k$. To adjust the carbon emissions of the fleet, a threshold value $H$ is considered to ensure that the average carbon emissions of the purchased vehicles is less than this value. %Additionally, each service region $i$ has its dedicated parking capacity for each car type $k$, which is represented by $C_i^k$, similar to \citet{Chang2017}. 
%This assumption is especially relevant when EVs are among the types of the cars purchased. Because, 
Since car sharing companies that utilize EVs can have charging stations at the parking locations of these vehicles, we consider parking stations dedicated to EVs to have charging facilities \citep{Chang2017, Brandstatter2020}. 
This further allows the company to start their operational daily planning with fully charged EVs which can remove the need for additional recharging throughout that day considering the longer driving ranges. %\esra{To this end, the parking space capacity of service region $i \in I$ for car type $k \in K$ is denoted by $C_i^k$.}

%We approach this problem with two different models. These models can be written in two-stage stochastic integer problem settings. In each settings, first stage problems include region opening and car allocation decisions. On the other hand, second stage problems include the movement of the cars in a spatial-temporal network. By using spatial-temporal network, car movements can be tracked since they are represented as flows in the network. For a car sharing company, it is important to see the time based movement of the cars, new cars can be added to the system, or some cars can be relocated to different region to satisfy predicted demand. 

%To maximize its profit, over $T$ time periods, a carshare company must decide number of regions that needs to be opened, denoted by $I$, distribute a given budget of $B$ vehicles among a set of opened regions. Let $K$ be the representation of different car type segments and $c_k$ be the price of purchasing a car type $k$. We assume that region $i \in I$ has a parking space capacity with respect to a car type $k \in K$, denoted by $C_i^k$, and $f_i$ is the fixed cost of the region, includes both opening cost of the region and parking space costs. 

We consider the operational level problem over $T$ time periods under demand uncertainty. These time periods correspond to the subperiods of a representative day for capturing the daily operations. 
To incorporate the operational level problem to the strategic level problem, the net profit obtained by this problem is scaled by multiplying it with $D$, where $D$ represents the number of operational days considered in the annual planning. 
The uncertainty is represented through demand scenarios, which are captured by set $W$, that are sampled from a given distribution. 
Each scenario $w \in W$ has a probability of occurrence $\pi_w$. 
We represent the customer demand through one-way trips and round-trips between service regions and time periods. In particular, for trips starting at period $t \in \{0, 1, \ldots, T-1\}$ and ending at period $s \in \{1, \ldots, T\}$, the parameter $d_{ijktsw}$ represents the demand for one-way trips from region $i \in I$ to region $j \in I \setminus \{i\}$ of car type $k \in K$ in scenario $w \in W$, and the parameter $d_{iktsw}$ represents the demand for round-trips for region $i \in I$ of car type $k \in K$ in scenario $w \in W$. 

To characterize the movement of the cars in the operational level problem, we construct a spatial-temporal network $G = (N,A)$ with a node set $N$ and an arc set $A$. Each node corresponds to a service region and time pair, which is denoted in the form of $n_{it}$ representing region $i \in I$ at period $t \in \{0, 1, \cdots, T\}$, where $t=0$ represents the status at the beginning of the operational planning. The directed arcs in this network indicate the movement of cars over time and space from one region to another from one time period to another. This network uses arcs of four different types as follows:
\begin{itemize}
\item One-way arcs in the form  $(n_{it},n_{js})$ correspond to the car flows of one-way trips from region $i$ to region $j$ from period $t$ to period $s$. The capacity of this arc in scenario $w$ for each car type $k$ depends on the demand amount $d_{ijktsw}$. 
\item Round-trip arcs in the form $(n_{it},n_{is})$ correspond to the car flows of round-trips for region $i$ from period $t$ to period $s$. The capacity of this arc in scenario $w$ for each car type $k$ depends on the demand amount $d_{iktsw}$. 
\item Relocation arcs in the form $(n_{it},n_{j,t+\zeta_{ij}})$ correspond to the car flows organized by the car sharing company to ensure rebalancing from region $i$ to region $j$ from period $t$ to $t+\zeta_{ij}$. Here, $\zeta_{ij}$ denotes the time that is needed to travel from region $i$ to $j$. 
These arcs are assumed to have sufficiently large capacity to ensure relocation operations. 
%These arcs have $\infty$ capacity and cost $c^{rel} \zeta_{ij}$ per trip.
\item Idle arcs in the form $(n_{it},n_{i,t+1})$ for $t = 1, \ldots, T-1$, correspond to the cars that are not used at region $i$ at the end of period $t$ after considering the flows on the relevant one-way arcs, round-trip arcs and relocation arcs. The capacities of these arcs for each car type $k$ is the corresponding parking capacity of the region $i$, which is $C_i^k$.  
\end{itemize}

\begin{table}[h]
\centering
\caption{Capacities and Unit Flow Revenues of Arc Types}
\label{tab:SpatialTemporalNetwork}
        \begin{tabular}{ccc}
        \hline
Arc Type                         & Capacity ($u_{akw}$) & Revenue ($r_{ak}$) \\ \hline
Idle Arc $a=(n_{it},n_{i,t+1})$ & $ C_i^k $  & 0                                       \\
One-Way Arc $a=(n_{it},n_{js})$              & $d_{ijktsw} $ & $r^{one}_k(s-t)$                               \\
Round-Trip Arc $a=(n_{it},n_{is})$           & $d_{iktsw}$  & $r^{two}_k(s-t)$                                 \\
Relocation Arc $a=(n_{it},n_{j,t+\zeta_{ij}})$           & $\infty $    & $-r^{rel}\zeta_{ij}$                    \\\hline
\end{tabular}
\end{table}

We denote the sets of one-way, round-trip, relocation and idle arcs by $A^{one}$, $A^{two}$, $A^{rel}$, $A^{idle}$, respectively, where $A$ represents the union of these four arc sets. Furthermore, to represent the arcs whose origin node is $n_{it}$, we define the set $\sigma^+(n_{it})$, and for the arcs whose destination node is $n_{it}$, we define the set $\sigma^-(n_{it})$. Table \ref{tab:SpatialTemporalNetwork} provides a summary of each arc type belonging to the spatial-temporal network with their arc capacities and unit flow revenues. 
For representing the arc capacities of each arc $a$ for car type $k$ in scenario $w$, we define the parameter $u_{akw}$. 
For representing the unit revenue of each arc $a$ for car type $k$, we define the parameter $r_{ak}$. 
In terms of the costs, one-way and round-trip arcs return profit by satisfying customer demand, whereas relocation arcs incur cost due to the resources used by the car sharing company to ensure rebalancing between different regions and time periods. 
Since car sharing companies generally adopt a time-based payment system that prices the trips based on their durations, the profit of one-way and round-trip arcs are computed based on the rental duration. Here, $r^{one}_k$ and $r^{two}_k$ represent the revenue of a car type $k$ per time unit for one-way trips and round-trips, respectively. Similarly, the cost of relocation arcs depend on the duration of the arc multiplied by the unit time cost of relocating cars, denoted by $r^{rel}$. 
For the idle arcs, company does not incur any additional cost, as the parking costs are included in the fixed cost of operating the service regions. 
%\textcolor{black}{Eger nodelara $k$ indexi de gelirse bu sekli bir $k$ icin gosteririz. Diger problem icin de sekil ekleyecek olsak uc boyutlu olmasi lazim, $k$'lar arasi baglantiyi gostermek icin. Ama bu durumda da flow variable'i $y$'de $k$ indexi gerekmiyor sanki, arc $k$ bilgisini de implicitly icereceginden?}    
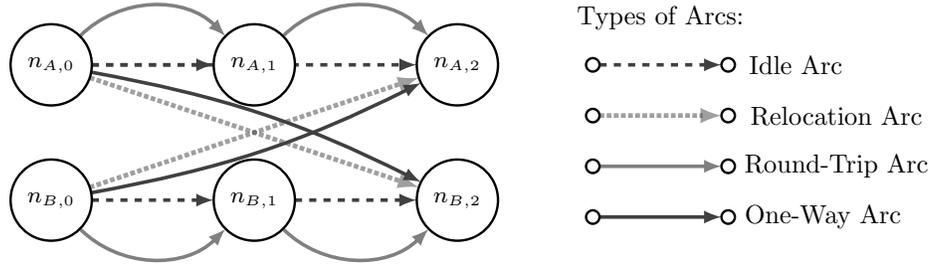
\begin{figure}[H]
\begin{center}
\caption{\centering A Spatial-Temporal Network Example} %\textcolor{black}{for a Given Car Type}
\label{fig:SpatialTemporalNetwork}
\scalebox{0.9}{
\begin{tikzpicture}
\Vertex[x=1,size=1.2,label=$n_{A, 0}$,fontscale=1.2,color=white]{A} \Vertex[x=4,size=1.2,label=$n_{A, 1}$,fontscale=1.2,color=white]{B} \Vertex[x=7,size=1.2,label=$n_{A, 2}$,fontscale=1.2,color=white]{C}
\Vertex[x=1,y=-2,size=1.2,label=$n_{B, 0}$,fontscale=1.2,color=white]{D}
\Vertex[x=4,y=-2,size=1.2,label=$n_{B, 1}$,fontscale=1.2,color=white]{E}
\Vertex[x=7,y=-2,size=1.2,label=$n_{B, 2}$,fontscale=1.2,color=white]{F}

%Idle Arcs
\Edge[Direct,style=dashed](A)(B)
\Edge[Direct,style=dashed](A)(C)
\Edge[Direct,style=dashed](D)(E)
\Edge[Direct,style=dashed](E)(F)
%Relocation
\Edge[Direct,style=densely dotted,opacity=0.5,lw=2](A)(F)
\Edge[Direct,style=densely dotted,opacity=0.5,lw=2](D)(C)
%Round-Trip
\Edge[Direct,bend=45,color=gray](A)(B)
\Edge[Direct,bend=45,color=gray](B)(C)
\Edge[Direct,bend=-45,color=gray](D)(E)
\Edge[Direct,bend=-45,color=gray](E)(F)
%One-Way
\Edge[Direct,opacity=1.2,bend=8](A)(F)
\Edge[Direct,opacity=1.2,bend=-8](D)(C)

%Legend
\Text[x=10,y=0.7]{Types of Arcs:}
\Text[x=12]{Idle Arc}
\Text[x=12.6,y=-0.75]{Relocation Arc}
\Text[x=12.6,y=-1.5]{Round-Trip Arc}
\Text[x=12.4,y=-2.25]{One-Way Arc}

\Vertex[x=9,size=0.01,color=white]{LA1}
\Vertex[x=11,size=0.01,color=white]{LA2}
\Vertex[x=9,y=-0.75,size=0.01,color=white]{LB1}
\Vertex[x=11,y=-0.75,size=0.01,color=white]{LB2}
\Vertex[x=9, y=-1.5,size=0.01,color=white]{LC1}
\Vertex[x=11, y=-1.5,size=0.01,color=white]{LC2}
\Vertex[x=9, y=-2.25,size=0.01,color=white]{LD1}
\Vertex[x=11, y=-2.25,size=0.01,color=white]{LD2}

\Edge[Direct,style=dashed](LA1)(LA2)
\Edge[Direct,style=densely dotted,opacity=0.5,lw=2](LB1)(LB2)
\Edge[Direct,color=gray](LC1)(LC2)
\Edge[Direct,opacity=1.2](LD1)(LD2)
\end{tikzpicture}}
\end{center}
\end{figure}
Figure \ref{fig:SpatialTemporalNetwork} provides a visual representation of the spatial-temporal network with two regions $A$ and $B$ over three time periods $\{0,1,2\}$ with sample one-way and round-trips. 
We note that each car type $k \in K$ corresponds to a different commodity that flows through our spatial-temporal network with 
its own set of one-way, round-trip, relocation and idle arcs. 
%This can be considered as a multicommodity flow representation as each arc type has a commodity corresponding to its car type. 
When substitution is integrated into this spatial-temporal network, then the demand of a customer for a specific car type can be satisfied by an another type of car, which implies that one-way and round-trip arcs can be shared across different commodities in the network. To address this issue, we define additional commodities as we introduce the operational problem with substitution in Section \ref{sec:SROP-S} by extending the flows in the arcs from the commodity of each car type to its multiple commodity types to explicitly capture the potential substitutions between different types of cars for satisfying the customer demand.

\subsection{Service Region and Operational Planning Problem (SROP)}
\label{sec:SROP}
%\textcolor{black}{BB: Buradan devam edecegim.}

To address the operational planning problem of the car sharing company, we define the integer decision variable $y_{akw}$, which represents the number of cars of type $k \in K$ in scenario $w \in W$ that are flowing on arc $a \in A$ over the spatial-temporal network. Combining above, we formulate the service region and operational planning (SROP) problem as follows:
\begin{subequations}
\label{eq:SROP-SingleStage}
\begin{alignat}{1}
\max \quad & - \sum_{i \in I} f_{i} z_{i} + D \sum_{w \in W} \sum_{a \in A}\sum_{k \in K} \pi_{w} r_{ak} y_{akw}   \label{eq:M1_1}   \\
\text{s.t.} \quad & x_{ik}  \leq  C_i^k z_{i} \quad \quad \forall i \in I, \> \forall k \in K, \label{eq:M1_2} \\
& \sum_{i \in I} \sum_{k \in K} c_{k} x_{ik} \leq B, \label{eq:M1_3} \\
& \sum_{i \in I}\sum_{k \in K} e_{k}  x_{ik} \leq H \sum_{i \in I}\sum_{k \in K} x_{ik} , \label{eq:M1_4}  \\
& \sum_{a \in \sigma^+(n_{it})} y_{akw} - \sum_{a \in \sigma^-(n_{it})} y_{akw} = \begin{cases} 
        x_{ik} & \text{if } t=0 \\
        0 & \text{if } t \in \{1,..., T-1\} \\
        -x_{ik} & \text{if } t=T
        \end{cases} 
         \quad \forall i \in I, \> \forall k \in K, \> \forall w \in W \label{eq:M1_6}  \\ 
& y_{akw} \leq u_{akw} z_{i} \quad \forall a=(n_{it},n_{is})  \in  A^{two}, \> \forall k \in K, \>\forall w \in W, \label{eq:M1_7}  \\
& y_{akw} \leq u_{akw} z_{i} \quad \forall a=(n_{it},n_{js}) \in A^{one}, \> \forall k \in K, \>\forall w \in W, \label{eq:M1_8} \\
& y_{akw} \leq u_{akw} z_{j} \quad \forall a=(n_{it},n_{js}) \in A^{one}, \> \forall k \in K, \>\forall w \in W, \label{eq:M1_9}  \\
&  y_{akw} \leq C_i^k z_{i} \quad \forall a=(n_{it},n_{i,t+1}) \in A^{idle},  \> \forall k \in K, \> \forall w \in W, \label{eq:M1_10} \\
&  z_{i} \in \{0,1\} \quad \quad \forall i \in I, \label{eq:M1_domain1}\\
&  x_{ik}  \in \mathbb{Z}^{+} \quad \quad \forall i \in I, \> \forall k \in K, \label{eq:M1_domain2}\\
& y_{akw} \in \mathbb{Z}^{+} \quad \quad \forall a \in A, \>  \forall k \in K, \> \forall w \in W.
\end{alignat}
\end{subequations}

The objective function \eqref{eq:M1_1} maximizes the total profit by considering the cost of operating service regions, expected revenue obtained from satisfying customer demand and expected cost of rebalancing cars between opened service regions. 
Constraint \eqref{eq:M1_2} ensures that the number of cars allocated to each open service region for each car type does not exceed their associated parking capacities. 
Constraint \eqref{eq:M1_3} corresponds to the total purchasing budget of the company in constructing its mix fleet of cars. 
Constraint \eqref{eq:M1_4} limits the total adjusted amount of $CO_2$ emissions by restricting the average emission amount over the purchased vehicles. In other words, the  average carbon emission of the fleet which is given by $\frac{\sum_{i \in I}\sum_{k \in K} e_{k}  x_{ik}}{\sum_{i \in I}\sum_{k \in K} x_{ik}}$ should not exceed the unit carbon allowance $H$ for each car the company owns. 

Constraints \eqref{eq:M1_6} are the flow balance constraints over the spatial-temporal network under every scenario $w \in W$. 
More specifically, at the beginning of the planning for the operational problem with $T$ time periods, each open service region $i \in I$ has $x_{ik}$ cars for each car type $k \in K$. During the intermediate time periods $t \in \{1, \cdots, T-1\}$, the number of cars leaving each node of the network is equal to the number of cars entering that node for every car type $k \in K$. 
At the last operational time period $T$, each open service region $i \in I$ has $x_{ik}$ cars for each car type $k \in K$, returning to their initial allocation. 
This consideration is done to have the same initial number of cars at each region for operating the car share system every $T$ periods. Moreover, as EVs are considered, allowing these vehicles to be fully charged at their parking spots at the end of $T$ time periods is necessary.%\esra{} 

%\textcolor{black}{BB: Paraphrase and add the following to make sure that this assumption is essential for charging of EVs: "We require that the locations of EVs in the last period are reset as their initial locations and all EVs are charged full, for the purpose of operating the carsharing system every T periods with the same initial deployment of fully charged EVs."}
 
Constraint \eqref{eq:M1_7} guarantees that the flow on the arcs corresponding to the round-trips should not exceed the capacity of their arcs over the open service regions. Similarly, constraints \eqref{eq:M1_8} and \eqref{eq:M1_9} are for the flows corresponding to the one-way trips by limiting the flow amounts on the relevant arcs when both origin and destination regions are opened.  
Constraint \eqref{eq:M1_10} considers the parking capacity of each open service region through the idle arcs for each car type. 
The remainder constraints ensure the integrality of the service region opening, fleet allocation and operational car flow decisions. 

Next in order, we reformulate the SROP problem \eqref{eq:SROP-SingleStage} as a two-stage stochastic program in \eqref{eq:SROP-TwoStage}, where the first-stage represents the strategic level planning problem and determines the service region opening $\mathbf{z}$ and fleet allocation decisions $\mathbf{x}$. Given these decisions, we define the second-stage problem for every scenario $w \in W$ as $\theta_{w}(\mathbf{z,x})$ optimizing the operational plan over the spatial-temporal network by maximizing the expected revenue minus the relocation costs. 
\begin{subequations}
\label{eq:SROP-TwoStage}
\begin{alignat}{1}
\max \quad &  - \sum_{i \in I} f_{i} z_{i} + D \sum_{w \in W} \pi_{w} \theta_{w}(\mathbf{z,x})  \\
\text{s.t.} \quad & \eqref{eq:M1_2} - \eqref{eq:M1_4}, \notag \\ %x_{ik}  \leq  C_i^k z_{i} \quad \forall i \in I, \> \forall k \in K,  \\
%& \sum_{i \in I} \sum_{k \in K} c_{k} x_{ik} \leq B,  \\
%& \sum_{i \in I}\sum_{k \in K} e_{k}  x_{ik} \leq H \sum_{i \in I}\sum_{k \in K} x_{ik},  \\
&  z_{i} \in \{0,1\}, \> x_{ik}  \in \mathbb{Z}^{+} \quad \forall i \in I, \> \forall k \in K, 
\end{alignat}
\end{subequations}
 where  for each scenario $w$;
 \begin{subequations}
 \label{eq:SROP-TwoStage-2nd}
\begin{alignat}{1}
\theta_{w}(\mathbf{z,x})  = & \max \sum_{a \in A} \sum_{k \in K} r_{ak} y_{akw} \label{eq:M1_5}\\
\text{s.t.} \quad & \sum_{a \in \sigma^+(n_{it})} y_{akw} - \sum_{a \in \sigma^-(n_{it})} y_{akw} = \begin{cases} 
        x_{ik} & \text{if } t=0 \\
        0 & \text{if } t \in \{1,..., T-1\} \\
        -x_{ik} & \text{if } t=T
        \end{cases} 
         \quad \forall i \in I, \> \forall k \in K,   \label{eq:M1_11} \\
        & y_{akw} \leq u_{akw} z_{i} \quad \forall a=(n_{it},n_{is})  \in  A^{two} ,\> \forall k \in K,  \label{eq:M1_12} \\ %, \> \forall i \in I silindi
        & y_{akw} \leq u_{akw} z_{i} \quad \forall a=(n_{it},n_{js}) \in A^{one} ,\> \forall k \in K, \label{eq:M1_13} \\ %, \> \forall i,j \in I, i \neq j silindi
        & y_{akw} \leq u_{akw} z_{j} \quad \forall a=(n_{it},n_{js}) \in A^{one} ,\> \forall k \in K, \label{eq:M1_14}  \\ %, \> \forall i,j \in I, i \neq j silindi
        &  y_{akw} \leq C_i^k z_{i} \quad \forall a=(n_{it},n_{i,t+1}) \in A^{idle}, \> \forall k \in K, \label{eq:M1_15} \\ %, \> \forall i \in I, \> \forall t \in \{1,..., T-1\} silindi
        & y_{akw} \in \mathbb{Z}^{+} \qquad \forall a \in A, \> \forall k \in K. 
\end{alignat}
\end{subequations}

This reformulation allows the second-stage problem to be solved efficiently by relaxing the integrality assumption on the flow variables $\mathbf{y}$ as follows. 

\begin{prop} \label{prop:TU-SROP}
Constraint matrix of the subproblem \eqref{eq:SROP-TwoStage-2nd} under any scenario $w \in W$ is totally unimodular.
\end{prop}

\proof{Proof:} 
Please see Appendix \ref{proof:TU-SROP}. 
\Halmos\endproof

\subsection{Service Region and Operational Planning Problem with Substitution (SROP-S)}
\label{sec:SROP-S}

A challenge for the SROP problem is that the demand for one-way trips and round-trips are lost if there are no available cars in the type requested by the customers. Nevertheless, in practice, car sharing companies can substitute between different car types by providing alternative options to the customers in case the specific car type requested by the customer is not available. 
This leads to flexibility in operations, higher customer satisfaction and better quality of service, which are desirable. 
To address this issue within planning, we introduce substitution to the SROP problem by allowing demand of one car type to be satisfied by another car type. 

%We define an extension for our SROP with using substitution between different car types. If the demand for a car type cannot be satisfied from that car type, and if there is an available car from other car types, substitution can be done. However, if there is no available car for a region and a time period, demand is lost. By using substitution among different car types, car sharing company can achieve higher customer satisfaction rates. 

To model this problem over the spatial-temporal network, we introduce additional  commodities in the form of a flow of a specific car type that is used for satisfying demand of the same or an alternative type of car. 
We define the set of commodities as $L$, where each commodity $l \in L$ can be represented as the flow of car type $k^1 \in K$ that is used for satisfying the demand of car type $k^2 \in K$. To this end, if $k^1$ is equal to $k^2$, then the demand of that car type is satisfied by itself, which is the case with no substitution. However, if $k^1$ is different than $k^2$, then substitution occurs by satisfying the demand of car type $k^2$ with a different type $k^1$. 
Subsequently, we propose two different commodity sets for each car type $k \in K$. 
First, let $F^k$ represent the set of commodities that use car type $k \in K$. Secondly, let $\hat{F}^k$ represent the set of commodities that can be used to satisfy the demand for car type $k \in K$. 
We note that for every car type $k \in K$, the set $F^k \cap \hat{F}^k$ returns a single commodity which represents the case without substitution where the demand of car type $k$ is satisfied by itself. Hence, the commodities in the sets $\cup_{k \in K}(F^k \cap \hat{F}^k)$ and $L \setminus \cup_{k \in K}(F^k \cap \hat{F}^k)$ represent the no substitution and substitution cases, respectively.
Additionally, the commodity set $L$ can be written as $L = \cup_{k \in K} F^k = \cup_{k \in K} \hat F^k$. 
To limit the substitution amount and potential customer dissatisfaction, we define a penalty parameter for substituted rentals as $p_l$ per time unit for every %$l \in \cup_{k \in K} (\hat{F}^k \setminus (F^k \cap \hat{F}^k))$ 
$l \in L \setminus \cup_{k \in K}(F^k \cap \hat{F}^k)$, which is discounted from rental price of the car. 
Thus, for each car type $k \in K$, we revise the revenue parameter defined in Table \ref{tab:SpatialTemporalNetwork} by introducing $r'_{al}$ for every $a \in A^{one} \cup A^{two}$ and $l \in L$, where $r'_{al} = r_{ak}$ for $l \in F^k \cap \hat{F}^k$, and $r'_{al} = r_{ak} - p_l$ for $l \in \hat{F}^k \setminus (F^k \cap \hat{F}^k)$. 
Since substitution is utilized for satisfying demand, the idle and relocation arcs are only defined over the commodities $l \in \cup_{k\in K} (F^k \cap \hat{F}^k)$, corresponding to the set of commodities representing the flows with no substitution, which further represents the commodities considered in the SROP problem. 
Thus, for each car type $k \in K$, for every $a \in A^{idle} \cup A^{rel}$, $r'_{al} = r_{ak}$ for $l \in F^k \cap \hat{F}^k$.

For illustrating the commodities and integration of substitution to the operational planning problem, we provide the following example setting with two car types. 

% edit etmeden onceki haline donmek istersek: end{document}in altinda
\begin{example}[Multicommodity Flows with Substitution]
%noindent{\bf Example (Multicommodity Flows with Substitution):}
We consider a car sharing company that classifies its cars into two types based on its fueling technology. 
More specifically, let $k = E$ represent electric vehicles and let $k = G$ represent gasoline vehicles. As we consider substitution, the demand of each vehicle type can be satisfied by the other vehicles. 
To this end, we define four commodities to represent these relationships:
\begin{itemize}
\item Commodity E-E: Electric cars used to satisfy the demand for electric cars
\item Commodity E-G: Electric cars used to satisfy the demand for gasoline cars
\item Commodity G-G: Gasoline cars used to satisfy the demand for gasoline cars
\item Commodity G-E: Gasoline cars used to satisfy the demand for electric cars
\end{itemize}

First, observe that commodities E-E and G-G correspond to the car flows with no substitution where the demand of each car type is satisfied by the same type of car. On the other hand, commodities E-G and G-E represent the substitution possibilities by allowing the demand of each car type to be satisfied by the other one. 
By using these four different commodities, we define $F^E = \{\text{E-E}, \text{E-G}\}$ for electric cars and $F^G = \{\text{G-G}, \text{G-E}\}$ for gasoline cars, since commodities E-E and E-G use electric cars, whereas commodities G-G and G-E use gasoline cars. 
On the other hand, we define $\hat{F}^E = \{\text{E-E}, \text{G-E}\}$ for electric cars and $\hat{F}^G = \{\text{E-G}, \text{G-G}\}$ for gasoline cars, to construct the sets of commodities that are used for satisfying the demand of each car type. 

Under this setting, the flows on the spatial-temporal network is in terms of these commodities, where one-way and round-trip arcs have car flows for each of the four commodities, and idle and relocation arcs include flows only for the commodities E-E and G-G. 
Table \ref{tab:SubstitutionExample} provides the objective function coefficients of each of the arc flows over the commodities in the operational problem by using the values of the revenue parameter $r'_{al}$ over every $a \in A$ and $l \in L$, whenever the corresponding commodity is defined on that arc. 
\end{example}
\begin{table}[H]
\centering
\caption{Objective Function Coefficients of the Multicommodity Arc Flows in the Operational Problem for the Example $(r'_{al})$}
\label{tab:SubstitutionExample}
\resizebox{\textwidth}{!}{%
\begin{tabular}{ccccc}
\hline
 & \multicolumn{2}{c}{Commodities Using Car Type E} & \multicolumn{2}{c}{Commodities Using Car Type G} \\ \cline{2-5} 
Arc types  & E-E                & E-G                   & G-G              & G-E                    \\ \hline
Idle $a=(n_{it},n_{i,t+1})$    & 0                & -                   & 0              & -                    \\
One-way $a=(n_{it},n_{js})$    & $r^{one}_{E}(s-t)$      & $(r^{one}_{E} - p_{E\text{-}G})(s-t)$ & $r^{one}_{G}(s-t)$    & $(r^{one}_{G} - p_{G\text{-}E})(s-t)$ \\
Round-trip $a=(n_{it},n_{is})$ & $r^{two}_{E}(s-t)$      & $(r^{two}_{E} - p_{E\text{-}G})(s-t)$ & $r^{two}_{G}(s-t)$    & $(r^{two}_{G} - p_{G\text{-}E})(s-t)$  \\
Relocation $a=(n_{it},n_{j,t+\zeta_{ij}})$ & $-c^{rel}_{E}\zeta_{ij}$ & -                   & $-c^{rel}_{G}\zeta_{ij}$ & -                    \\ \hline
\end{tabular}%
}
\end{table}

%To integrate substitution and the corresponding multicommodities to the 
For constructing the operational problem with substitution, we introduce the integer variable $y_{alw}$, for $a \in A^{one} \cup A^{two},\, l \in L$ and  $a \in A^{idle} \cup A^{rel}, \, l \in \cup_{k \in K} F^{k} \cap \hat{F}^{k}$, which indicates the number of cars of commodity $l$ in scenario $w \in W$ that are flowing on arc $a$ over the spatial-temporal network. 
We formulate the resulting Service Region and Operational Planning with Substitution (SROP-S) problem as a two-stage stochastic program in \eqref{eq:SROP-S-TwoStage}. 
Similar to the SROP problem, the first-stage problem optimizes the service region opening $\mathbf{z}$ and fleet allocation decisions $\mathbf{x}$. Given these decisions, we revise the second-stage problem for every scenario $w \in W$ by defining $\bar \theta_{w}(\mathbf{z,x})$ that allows substitution between different car types. 
\begin{subequations}
\label{eq:SROP-S-TwoStage}
\begin{alignat}{1}
\max \quad &  - \sum_{i \in I} f_{i} z_{i} + D \sum_{w \in W} \pi_{w} \bar \theta_{w}(\mathbf{z,x})  \\
\text{s.t.} \quad & \eqref{eq:M1_2} - \eqref{eq:M1_4}, \notag \\ %x_{ik}  \leq  C_i^k z_{i} \quad \forall i \in I, \> \forall k \in K,  \\
%& \sum_{i \in I} \sum_{k \in K} c_{k} x_{ik} \leq B,  \\
%& \sum_{i \in I}\sum_{k \in K} e_{k}  x_{ik} \leq H \sum_{i \in I}\sum_{k \in K} x_{ik},  \\
&  z_{i} \in \{0,1\}, \> x_{ik}  \in \mathbb{Z}^{+} \quad \forall i \in I, \> \forall k \in K,
\end{alignat}
\end{subequations}
 where for each scenario $w$;
 \begin{subequations}
 \label{eq:SROP-S-TwoStage-2nd}
\begin{alignat}{1}
\bar \theta_{w}(\mathbf{z,x})\quad &  = \max \sum_{a \in A^{one} \cup A^{two}} \sum_{l \in L} r'_{al} y_{alw} +
\sum_{a \in A^{rel}} \sum_{k \in K} \sum_{l \in F^{k} \cap \hat{F}^{k}} r'_{al} y_{alw} \label{eq:M2_obj} \\ 
%ilk sum su ekideydi:\max \sum_{a \in A^{one} \cup A^{two}} \sum_{k \in K} \sum_{l \in \hat{F}^k} r'_{al} y_{alw}
\text{s.t.} \quad & \sum_{l \in F^k} \sum_{a \in \sigma^+(n_{it})} y_{alw} - \sum_{l \in F^k} \sum_{a \in \sigma^-(n_{it})} y_{alw} = \begin{cases} 
        x_{ik} & \text{if } t=0 \\
        0 & \text{if } t \in \{1,..., T-1\} \\
        -x_{ik} & \text{if } t=T
        \end{cases} 
        \quad \forall i \in I, \> \forall k \in K,  \label{eq:M2_1} \\
        & \sum_{l \in \hat{F}^k} y_{alw} \leq u_{akw} z_{i} \quad \forall a=(n_{it},n_{is}) \in  A^{two} ,\> \forall k \in K, \label{eq:M2_2} \\ %, \> \forall i \in I  silindi 
        & \sum_{l \in \hat{F}^k} y_{alw} \leq u_{akw} z_{i} \quad \forall a=(n_{it},n_{js}) \in A^{one} , \> \forall k \in K,  \label{eq:M2_3}\\ %\> \forall i,j \in I, i \neq j, silindi
        &  \sum_{l \in \hat{F}^k} y_{alw} \leq u_{akw} z_{j} \quad \forall a=(n_{it},n_{js}) \in A^{one} , \> \forall k \in K, \label{eq:M2_4} \\ %\> \forall i,j \in I, i \neq j silindi
        &   y_{alw} \leq C_i^k z_{i} \quad \forall a=(n_{it},n_{i,t+1}) \in A^{idle},  \> \forall l \in F^{k} \cap \hat{F}^{k}, \> \forall k \in K,  \label{eq:M2_5} \\ %\> \forall i \in I, \> t \in \{1,..., T-1\}, silindi
        & y_{alw} \in \mathbb{Z}^{+} \quad \forall a \in A^{one} \cup A^{two}, \>  l \in L, \\ %l \in \hat{F}^k, \> \forall k \in K
        & y_{alw} \in \mathbb{Z}^{+} \quad \forall a \in A^{idle} \cup A^{rel}, \>  l \in F^{k} \cap \hat{F}^{k}, \> \forall k \in K.
\end{alignat}
\end{subequations}

Under each demand realization of scenario $w \in W$, objective function of the second-stage problem \eqref{eq:M2_obj} maximizes the revenue obtained from one-way and round-trips minus the relocation cost and penalty of substitution. 
Constraints \eqref{eq:M2_1} correspond to the flow balance constraints. Different than the second-stage problem of the SROP problem, for each car type $k \in K$, the flows over the multicommodities are summed over $F^k$, the set of commodities using car type $k$, to ensure the flow balance. Constraints \eqref{eq:M2_2} - \eqref{eq:M2_4} are analogous to constraints \eqref{eq:M1_12} - \eqref{eq:M1_14} in representing round-trips and one-way trips, whereas the flows for each car type $k \in K$ are summed over $\hat F^k$ to consider the set of car flows satisfying the demand of car type $k$. 
Lastly, constraint \eqref{eq:M2_5} ensures that for each car type $k \in K$, the idle cars at the end of each period are within the parking limit of the corresponding region and flow  through its regular commodity as defined through $F^{k} \cap \hat{F}^{k}$. 

%\subsection{TU Proof}
%\textcolor{blue}{EK: constraint 10u ana modeldeki idle arc icin olan kisit (5f) yerine kullanabiliriz. Bu arada sanki su anki 2nd model dualimizin idle arc icin olan kisitinda bir gariplik var gibi. O kisitta k'=k olmak zorunda gibi, baska turlusu mumkun degil sanki. Kontrol eder misin Sinan? k'=k oldugunda da zaten constraint 10un dualine denk geliyor. idle arc icin yeni dual kisitimiz:
%\begin{equation*}
%    \beta_{itk} - \beta_{i,t+1,k} + \lambda_{ak} \geq 0 \quad  a = (n_{it},n_{i,t+1}) \in A^{idle} \quad  l_k \in L, k \in K
%\end{equation*}

%where $l_k$ represents the commodity that uses car type $k$ for satisfying the demand for car type $k$.
%}

%\textcolor{blue}{where $l_k$ represents the commodity that uses car type $k$ fot satisfying the demand for car type $k$. In our example, $l_1 = 1$ and $l_2 = 3$.

%Eski halinde L kullanimisim, ama L set of commodities icin kullanilmis: Assume that $L$ is the subset of commodities where the demand is satisfied by the desired car type. For example, in our case $L$ is given by $\{1,3\}$. }  

\begin{prop} \label{prop:TU-SROP-S}
Constraint matrix of the subproblem \eqref{eq:SROP-S-TwoStage-2nd} under any scenario $w \in W$ is totally unimodular.
\end{prop}

\proof{Proof:} 
Please see Appendix \ref{proof:TU-SROP-S}. 
\Halmos\endproof

We benefit from the totally unimodularity of the constraint matrix of this second-stage problem as we design our exact solution algorithm to solve the challenging SROP-S problem. 

\begin{remark}
We note that the SROP-S problem is a generalization of the SROP problem. More specifically, when $F^k = \hat{F}^k$ for every car type $k \in K$, then each of these sets only include one commodity flow corresponding to the case when the demand of car type $k$ is satisfied by the car type $k$ itself. This reduces the SROP-S problem to the SROP problem as substitution is disregarded. 
Thus, by adjusting the car types allowed for substitution, the car sharing company can determine its flexibility level through leveraging the SROP-S problem. 
%This result also applied to the proof of Proposition EKLE. 
\end{remark}

%\textcolor{black}{In problem statement or the section with substitution, we can add a note that states that customers are assumed to accept the substitute trips due to the incentives/reduced pricing suggested by the car sharing company!}

%\textcolor{black}{Should we add a remark/note regarding relocation capacity as well if it does not violate TU property?}

\section{Solution Algorithm} \label{sec:solution}
The two-stage stochastic programming models presented in the previous section are large-scale mixed integer programs (MIP) due to the potentially large number of scenario-dependent decision variables and constraints. The most popular procedure for solving this type of large-scale stochastic optimization models is to use the L-Shaped method \citep{vanslyke1969} which is actually the application of Benders decomposition \citep{benders} for solving two-stage stochastic linear programs (LP). The main advantage of Benders decomposition for solving two-stage stochastic programs is that the second-stage problem decomposes for each scenario when the first-stage variables are fixed. By benefiting from our results in Propositions \ref{prop:TU-SROP} and \ref{prop:TU-SROP-S}, we make use of this decomposing structure of our stochastic MIPs to develop an efficient branch-and-cut algorithm with computational enhancements.

%MILP models can be very hard to solve when problem size gets larger and more variables and constraints are involved. Thus, we use Benders Decomposition to improve the solution times. Benders decomposition is a useful algorithm to solve MILP models with large number of variables and constraints. The algorithm works with the problems that can be divided into two stages. The first stage decisions are the master problem and the second stage decisions can be considered as separate sub-problems. The first attempt at first stage problem is made without any knowledge about optimality and feasibility of second stage problems. After the first attempt, master problem is fed by information from sub-problems. Then, the master problem is solved accordingly. This process provides lower and upper bound on the optimal value objective value iteratively. These iterations ends when the gap between lower and upper bound reachs zero, or there is no finite optimal solution. We showed that SROP and SROP-S can be written as two-stage models. Then, we can define an auxiliary variable $Q_w$ represents the value of the sub-problem $w$ at optimality.     

In the subsequent sections, we present our Benders decomposition algorithm for SROP and SROP-S, which is an iterative process by solving a relaxed master problem (RMP) and generating cuts to be added for this problem accordingly. For both SROP and SROP-S, we keep the first stage decision variables $\mathbf{z}$ and $\mathbf{x}$ in the Benders RMP, but as their second-stage problems are different, the subproblems solved and the cuts generated are different as explained below. 

\subsection{Benders Decomposition for SROP} \label{sec:bendersSROP}
In our decomposition algorithm for SROP, in addition to the first stage decision variables $\mathbf{z}$ and $\mathbf{x}$ that will be considered in RMP we define auxiliary decision variables $Q_w$ to approximate the optimal value of the subproblem for scenario $w \in W$. Accordingly, at an intermediate step of the algorithm, we solve the following RMP:  
\begin{subequations}
\begin{alignat}{1}
\max \quad &  - \sum_{i \in I} f_{i} z_{i} + D \sum_{w \in W} \pi_{w} Q_{w}\\  
\text{s.t.} \quad & \eqref{eq:M1_2} - \eqref{eq:M1_4}, \eqref{eq:M1_domain1}, \eqref{eq:M1_domain2} \notag \\%x_{ik}  \leq  C_i^k z_{i} \qquad \forall i \in I, \> \forall k \in K,\\
%& \displaystyle\sum_{i \in I}\sum_{k \in K} e_{k}  x_{ik} \leq \mathcal{H}\displaystyle \sum_{i \in I}\sum_{k \in K} x_{ik},\\
%&\displaystyle\sum_{i \in I}\sum_{k \in K} c_{k} x_{ik} \leq B, \\
&O(\mathbf{z,x,Q}) \geq 0, \\
& Q_w \leq U_w \quad \forall w \in W,
\end{alignat}
\end{subequations}
where $O(\mathbf{z,x,Q}) \geq 0$ represents the optimality cuts generated and added to RMP until that iteration, and $U_w$ represents an upper bound for the net profit that can be obtained under scenario $w \in W$. Since keeping the cars idle in their service regions in all periods through the idle arcs is a feasible solution for the subproblem \eqref{eq:SROP-TwoStage-2nd} under any $\mathbf{z}$ and $\mathbf{x}$ solution, Benders feasibility cuts are not needed in our algorithm. $U_w$ can be easily determined as the total profit of satisfying all demand under scenario $w \in W$ by ignoring the decisions of the problem.  
%\textcolor{black}{Let us discuss adding upper bound on $Q_w$ values, in that case we might not need to add the nonnegativity constraint on $Q_w$.}

Given an optimal solution $(\hat{\mathbf{z}},\hat{\mathbf{x}},\hat{\mathbf{Q}})$ for RMP, we solve the dual of the second-stage problem \eqref{eq:SROP-TwoStage-2nd} as our subproblem for each scenario $w \in W$ in \eqref{eq:SROP-Benders-dualsub}. To this end, dual variables $(\boldsymbol{\beta}, \boldsymbol{\alpha}, \boldsymbol{\gamma^{1}}, \boldsymbol{\gamma^{2}}, \boldsymbol{\lambda})$ are associated with the constraints \eqref{eq:M1_11}, \eqref{eq:M1_12}, \eqref{eq:M1_13}, \eqref{eq:M1_14}, \eqref{eq:M1_15}, respectively. 
\begin{subequations}
 \label{eq:SROP-Benders-dualsub}
\begin{alignat}{1}
\begin{split}  
\displaystyle \min \quad&  \sum_{i\in I} \sum_{k\in K}  \hat{x}_{ik} (\beta_{i0kw} - \beta_{iTkw})  + \sum_{i \in I}\sum_{a \in  A^{two}(i)}\sum_{k \in K} u_{akw} \hat{z}_{i} \alpha_{akw}  + 
\sum_{i \in I}\sum_{a \in  A^{one}(i^{+})} 
\sum_{k \in K} u_{akw} \hat{z}_{i} \gamma^{1}_{akw} \\ 
& +  \sum_{i \in I}\sum_{a \in  A^{one}(i^{-})}\sum_{k \in K} u_{akw} \hat{z}_{i} \gamma^{2}_{akw} +  
     \sum_{i \in I} \sum_{k \in K} C_i^k \hat{z}_{i} (\sum_{a \in A^{idle}(i)} \lambda_{akw})
\end{split}\\
\text{s.t.}\quad
& \beta_{itkw} - \beta_{jskw} + \gamma^{1}_{akw}+ \gamma^{2}_{akw}  \geq r_{ak} \quad  a = (n_{it},n_{js}) \in A^{one}, \>  k \in K, \\
& \beta_{itkw} - \beta_{iskw} +  \alpha_{akw}  \geq r_{ak} \quad  \forall a = (n_{it},n_{is}) \in A^{two}, \> k \in K, \\
& \beta_{itkw} - \beta_{j,t+\zeta_{ij},kw} \geq r_{ak} \quad  \forall a = (n_{it},n_{j,t+\zeta_{ij}}) \in A^{rel}, \> k \in K, \\
& \beta_{itkw} - \beta_{i,t+1,kw} + \lambda_{akw} \geq 0 \quad \forall a = (n_{it},n_{i,t+1}) \in A^{idle}, \> k \in K, \\
&\gamma^{1}_{akw}, \gamma^{2}_{akw} \geq 0 \quad \forall a \in  A^{one}, \> k \in K, \\
& \alpha_{akw} \geq 0 \quad \forall a \in  A^{two}, \> k \in K, \\
& \lambda_{akw} \geq 0 \quad \forall a \in A^{idle}, \> k \in K,
\end{alignat} 
\end{subequations}
where $A^{two}(i)$, $A^{one}(i^+)$, $A^{one}(i^-)$, $A^{idle}(i)$ represent the subsets of arcs that are defined for region $i \in I$. More specifically, $A^{two}(i) = \{a = (n_{it},n_{is}) \in A^{two}: t,s \in \{0,\ldots, T\}, t < s \}$ represents the set of all round-trip arcs incident to nodes defined for region $i \in I$,  $A^{one}(i^+) = \{a = (n_{it},n_{js}) \in A^{one}: t,s \in \{0,\ldots, T\}, t < s, \; j \in  I\}$ and $A^{one}(i^-) = \{a = (n_{jt},n_{is}) \in A^{one}: t,s \in \{0,\ldots, T\}, t < s, \; j \in  I\}$ represent the set of all one-way arcs with the origin and destination, respectively, of region $i \in I$, and $A^{idle}(i) = \{a = (n_{it},n_{i,t+1}) \in A^{idle}: t \in \{0,\ldots, T-1\}\}$ gives the set of all idle arcs defined for region $i \in I$. 

At each iteration of the algorithm, we first solve the RMP and get its optimal solution $(\hat{\mathbf{z}},\hat{\mathbf{x}},\hat{\mathbf{Q}})$. Given this solution, we solve the subproblem \eqref{eq:SROP-Benders-dualsub} and obtain its optimal solution as $(\hat{\boldsymbol{\beta}},\hat{\boldsymbol{\alpha}}, \hat{\boldsymbol{\gamma^{1}}}, \hat{\boldsymbol{\gamma^{2}}}, \hat{\boldsymbol{\lambda}})$ for every scenario $w \in W$. %be the optimal solution of the dual subproblem for scenario $w \in W$ given by \eqref{eq:SROP-Benders-dualsub}. 
If its objective function value is overestimated in the current optimal solution %$(\mathbf{z},\mathbf{x}, \mathbf{Q})$ 
of RMP, i.e. if the optimal value of the subproblem for $w$ is smaller than $\hat{Q}_w$, then we add the following optimality cut to RMP for scenario $w \in W$:
%Benders decomposition solves a relaxed master problem to obtain a canditate optimal solution $(z^{*},x^{*},\theta^{*})$ and use this solution to solve dual-problems to optimize $\theta_{w}(\mathbf{z^{*},x^{*}})$. If a sub-problem has an optimal solution that has lower value than the optimal value given by the master problem, then an optimality cut is added to the relaxed master problem to solve the relaxed master problem again. The optimality cut for this problem can be written as:
\begin{alignat}{1}
\begin{split}
& \displaystyle  \sum_{i\in I} \sum_{k\in K} x_{ik} (\hat{\beta}_{i0kw} - \hat{\beta}_{iTkw})  + \sum_{i \in I}\sum_{a \in  A^{two}(i)}\sum_{k \in K} u_{akw} z_{i} \hat{\alpha}_{akw}  +  \sum_{i \in I}\sum_{a \in  A^{one}(i^{+})}\sum_{k \in K} u_{akw} z_{i} \hat{\gamma}^{1}_{akw} \\
&+ \sum_{i \in I}\sum_{a \in  A^{one}(i^{-})}\sum_{k \in K} u_{akw} z_{i} \hat{\gamma}^{2}_{akw} +  
     \sum_{i \in I} \sum_{k \in K} C_i^k z_{i} (\sum_{a \in A^{idle}(i)} \hat{\lambda}_{akw}) - Q_{w} \geq 0
\end{split}
\end{alignat}
%\textcolor{black}{If no cut can be added for any scenario $w \in W$, then the algorithm terminates with the optimal solution, which is the solution obtained by the RMP.}

\subsection{Benders Decomposition for SROP-S} \label{sec:bendersSROP-S}
Similar to the previous section, to solve SROP-S with Benders decomposition we first introduce auxiliary decision variables $\bar{Q}_w$ to approximate the optimal value of the subproblem for scenario $w \in W$. Hence, the following RMP will be solved at an intermediate iteration of the algorithm for SROP-S:
\begin{subequations}
\begin{alignat}{1}
\max \quad &  - \sum_{i \in I} f_{i} z_{i} + D \sum_{w \in W} \pi_{w} \bar{Q}_{w}\\  
\text{s.t.} \quad & \eqref{eq:M1_2} - \eqref{eq:M1_4}, \eqref{eq:M1_domain1}, \eqref{eq:M1_domain2} \notag \\ %x_{ik}  \leq  C_i^k z_{i} \qquad \forall i \in I, \> \forall k \in K,\\
%& \displaystyle\sum_{i \in I}\sum_{k \in K} e_{k}  x_{ik} \leq \mathcal{H} \sum_{i \in I} \sum_{k \in K} x_{ik},\\
%&\displaystyle\sum_{i \in I}\sum_{k \in K} c_{k} x_{ik} \leq B, \\
&\Bar{\mathcal{O}}(\mathbf{z,x,\bar{Q}}) \geq 0, \\
& \bar{Q}_w \leq \bar{U}_w \quad w \in W, 
\end{alignat}
\end{subequations}
where $\Bar{\mathcal{O}}(\mathbf{z,x,\bar{Q}}) \geq 0$ includes optimality cuts generated until that iteration. Since SROP-S is a generalization of SROP, the feasible solution of keeping all vehicles idle in their initial regions for all periods is also a feasible solution for the subproblem of SROP-S. Hence, we do not need to consider feasibility cuts in our algorithm. 
The upper bound $\bar{U}_w$ can be computed in a similar manner as in SROP model. 

Given an optimal solution $(\hat{\mathbf{z}},\hat{\mathbf{x}}, \hat{\mathbf{\bar{Q}}})$ for RMP, we solve the dual of the second-stage problem \eqref{eq:SROP-S-TwoStage-2nd} as our subproblem for each scenario $w \in W$ by relating the dual variables $(\boldsymbol{\beta}, \boldsymbol{\alpha}, \boldsymbol{\gamma^{1}}, \boldsymbol{\gamma^{2}}, \boldsymbol{\lambda})$ with the constraints \eqref{eq:M2_1}-\eqref{eq:M2_5}, respectively. \textcolor{black}{Please see Appendix \ref{sec:DualSecondStageSROPS} model \eqref{eq:SROP-S-Benders-dualsub} for the dual of the second-stage problem.}  
%edit edilmeden onceki hali enddocument'den sonra

Let $(\hat{\boldsymbol{\beta}}, \hat{\boldsymbol{\alpha}}, \hat{\boldsymbol{\gamma^{1}}}, \hat{\boldsymbol{\gamma^{2}}}, \hat{\boldsymbol{\lambda}})$ be the optimal solution of the dual subproblem for scenario $w \in W$ given by \eqref{eq:SROP-S-Benders-dualsub}. If its objective function value is overestimated in the current optimal solution $(\hat{\mathbf{z}},\hat{\mathbf{x}}, \hat{\mathbf{\bar{Q}}})$ of RMP, i.e. if the optimal value of the subproblem for $w$ is smaller than $\hat{\bar{Q}}_w$, then we add the following optimality cut to RMP for scenario $w \in W$:
\begin{alignat}{1}
\begin{split}
& \displaystyle  \sum_{i\in I} \sum_{k\in K} x_{ik} (\hat{\beta}_{i0kw} - \hat{\beta}_{iTkw})  + \sum_{i \in I}\sum_{a \in  A^{two}(i)}\sum_{k \in K} u_{akw} z_{i} \hat{\alpha}_{akw}  +  \sum_{i \in I}\sum_{a \in  A^{one}(i^{+})}\sum_{k \in K} u_{akw} z_{i} \hat{\gamma}^{1}_{akw} \\
&+ \sum_{i \in I}\sum_{a \in  A^{one}(i^{-})}\sum_{k \in K} u_{akw} z_{i} \hat{\gamma}^{2}_{akw} +  
     \sum_{i \in I} \sum_{k \in K} C_i^k z_{i} (\sum_{a \in A^{idle}(i)} \hat{\lambda}_{akw}) - \bar{Q}_{w} \geq 0
\end{split}
\end{alignat}
%\textcolor{black}{BB: We should also mention branch-and-cut implementation, average scenario cuts, warm start, how we improve lower and upper bounding procedures etc. as we present our computational enhancements.}

\subsection{Computational Enhancements}\label{sec:imp&enhanc} %
To enhance our solution algorithm, we implement the Bender's decomposition algorithms by building a single search tree for the RMP. This is accomplished by employing the lazy constraint callback feature of the off-the-shelf solver. When a new incumbent solution is found in the branch-and-cut tree, the lazy constraint callback is invoked to solve the subproblems for each scenario for the current incumbent solution. If the optimal value of the subproblem for $w \in W$ is smaller than $Q_w$ (or $\bar{Q}_w$), then an optimality cut is added to the RMP. If no cut is violated for any scenario, then the current solution is considered as a new incumbent solution. Following the literature, we use the multi-cut version of the algorithm where an optimality cut is added separately for each subproblem \citep{rahmaniani2017}. 

The naive implementation of the Benders or L-Shaped algorithm mostly
does not perform well because of the information lost in RMP due to the decision variables and constraints removed to the subproblems. \cite{rahmaniani2018} summarizes, tests and compares different acceleration strategies that are proposed and used in the literature to overcome this drawback of the algorithm. In our computational experiments, we test the following acceleration strategies for our branch-and-cut algorithms for SROP and SROP-S:
\begin{enumerate}[(i)]
    \item \textbf{Initial solution:} As a warm-start strategy, we provide initial solutions for the problems by solving relatively small-scale problems with a randomly selected single scenario. In other words, we solve the MIP formulations of the problems by assuming that we have a single scenario, and we use the first stage solution obtained from these MIPs as an initial solution for our RMP.  % average scenario that represents expected values of the random parameters: bunu yapmadik sanirim. 
    \item \textbf{Initial Benders cuts:} Given the initial solution found in the previous item, we construct the corresponding Benders optimality cuts for each subproblem, and add these cuts to RMPs as initial cuts. These cuts enable the algorithm to have better (smaller) upper bounds in the earlier iterations.
    \item \textbf{Improving the LP relaxation: } Lazy constraint callback and user cut callback are two important features of the off-the-shelf solvers that add cuts to the problems at different times. Lazy constraint callback is called when an integer candidate solution is found for the RMP while the user cut callback is called at any fractional solution. Since the Benders' cuts are valid inequalities for RMPs we add them also at fractional solutions using the user cut callback feature of the off-the-shelf solver. Following the literature, instead of adding user cuts at any node of the branch-and-cut tree, we add them only in the root node until the improvement in the relative optimality gap is very small. In this way, we improve the LP relaxations of RMPs at the root node of the branch-and-cut tree as much as possible.
\end{enumerate}

\section{Computational Study}
\label{sec:ComputationalStudy}

Our computational study demonstrates the value of the proposed models in terms of the service region, fleet design and operational planning decisions with the integration of substitution, along with the computational efficiency of the solution approach. In Section \ref{sec:insights}, from a modeling point of view, we highlight the impact of our modeling approach and the key problem parameters on optimal solutions from various perspectives with managerial insights. In Section \ref{sec:performance}, we investigate the effectiveness of our decomposition-based algorithms compared to the solution of MIPs by a commercial solver. Before proceeding to the results, we first describe the generation of problem instances used in our experiments in Section \ref{sec: data}.

\subsection{Experimental Setup} \label{sec: data}
We generate problem instances that represent the real case as much as possible by using the parameter settings explained in \cite{Luetal2017} which is based on the data set of Zipcar in the Boston-Cambridge, Massachusetts area. Different than this study, we consider a car sharing company that aims to determine service regions and build its fleet from a mixture of vehicle types. Specifically, we focus on a company building a fleet from two different car types $K=\{E,G\}$ where $E$ and $G$ represent the electrical and gasoline cars, respectively. We search for the prices and carbon emissions of different car models available in the market (see,\cite[e.g.][]{toyotaus, carbon-offsets, epa} ), and accordingly set $c = [34K, 27K]$ and $e = [0,0.75]$. We assume that the unitary carbon allowance is $H = 0.5$ unless otherwise is stated. Note that due to the carbon emission constraints \eqref{eq:M1_4}, $H=0.5$ implies that for every two gasoline cars purchased one electric car should be also purchased. We investigate the effect of the value of $H$ in the last part of Section \ref{sec:insights}. 
%These nine different region have 6, 9, 7, 6, 8, 9, 8, 9, and 6 parking space per car type. 
 We consider an area that is divided into 9 equal possible service regions, i.e. $|I|=9$. We generate the parking space capacities of each region for each car type $C_{i}^k$ randomly from $U[6,9]$ for $i \in I$ and $k \in K$. Similar to \cite{Luetal2017}, we set the annual unit parking space cost to $\$3500$ for each gasoline car in each region. For the electric car parking spaces, due to the charging system operations, we set a higher parking space cost which is $\$4000$. We assume that the fixed rental cost for each service region is $\$300K$. To determine the total fixed cost $f_i$ for opening service region $i \in I$ and locating the parking spaces, we add the total cost of the parking spaces to the rental cost of the region. For instance, for region 1, $f_1 = 300K + 6 \times (3500 +4000) =  345K$. \textcolor{black}{Please see Table \ref{tab:data1} in Appendix \ref{sec:ServiceRegionCapacities} for the parking space capacities and fixed costs for these regions.} % The parking space capacities and fixed costs for these regions can be seen in Table \ref{tab:data1}. 
  
 %In a real life application, the sizes and the rental (fixed) prices of the service regions might depend on many things. For instance, it should be more expensive to find a sufficiently large area for parking spaces in the center of a metropolitan area compared to its neighborhood \esra{(datamiz boyle degil aslinda, bu cumleyi kaldirabiliriz?)}. 
%Fixed costs are $\$345K$, $\$367.5K$, $\$352.5K$, $\$345K$, $\$360K$, $\$367.5K$, $\$360K$, $\$367.5K$, and $\$345K$ respectively.     

We consider a daily operational plan by dividing 24 hours into $|T|=12$ periods. This can be done in different ways. For instance, one can simply divide 24 hours into 12 equal length periods. Alternatively, one can represent the rush-hours by using more periods of shorter length while merging several hours into one period during the off-hours. We use the first one by assuming that each period has the length of 2 hours. Since we make a daily plan for the trips, we set $D=365$ as a normalization factor for combining the yearly fixed costs with the daily profits.  

We determine the relocation times (in terms of periods) between the service regions based on the distance between them such that $\zeta_{ij} = 1$ if regions $i$ and $j$ are neighbors, and $\zeta_{ij} = 2$ otherwise. In addition to the fact that it takes more time to travel between two non-neighbor regions, we consider this setting to discourage frequent relocation actions between two regions that are far away from each other. 

We consider the duration of the trips while generating the demand for the trips between regions. In other words, we assume that the probability of observing a demand for a long duration trip is very small compared to the shorter duration trips inspiring from the study of \cite{strohle2019} where it is empirically shown that majority of the trips are short in terms of duration and distance. More specifically, we assume that the demand for a one-way or round-trip with a duration less than or equal to $\frac{|T|}{3}$ is 0, 1, and 2 with the probabilities 0.8, 0.15, 0.05, respectively. If the duration of a trip is larger than $\frac{|T|}{3}$, since it is less likely to occur, the demand is $0$ with probability $0.8$, and $1$ with probability $0.2$. We use this demand distribution for generating all demand scenarios between all pairs of regions, and assume that the probability of realizing each scenario is the same, i.e. $ \pi_w = \frac{1}{|W|}$.  

Following the settings of \cite{Luetal2017}, we set the revenue for one-way and round-trips as $r_k^{one}=\$12$ and $r_k^{two}=\$7.75$ per hour for each car type $k$. We also assume that the relocation cost is $c^{rel}=\$8$ per hour. In the model SROP-S, we set the penalty cost for substitution to $p=2$ unless otherwise is stated. 
For illustrative purposes, in the subsequent section, we consider this problem instance with $|W|=100$ scenarios and analyse the results in detail from different perspectives.

\ignore{
{\color{red}BB: Bu sectionda her bir nodeun neye denk geldigini, kac gunluk planlama yaptigimizi, bir gunu time periyotlara nasil boldugumuzu vb. kisimlari da aciklayalim. Ayrica, aslinda second-stage problemi objective'de 365 ile carpiyorduk degil mi? O yuzden o gun sayisini da bir parametre tanimlayarak yukaridaki modellere ekleyelim. Parametrenin degerinin 365 oldugunu da burada soyleriz.}

{\color{red}BB: Computational study setting ve parametre degerlerinin ne alindigi, ilgili makalelere referans verilerek buraya eklenmeli. }

{\color{red}BB: Carbon emission parametresi ile hic oynamadik. Benzer sekilde first-stagede bulunan butce kisidi ile ilgili ama iki farkli settingde de run almis olduk. Acaba bununla ilgili de bazi runlar alsak mi? Su anki carbon parametre degeri altinda EV-Regular araba dagilimini nasil degerlendiriyoruz? Bu parametre artar azalirsa ne olur? Ornegin daha cok EV almaya zorlarsak, regular talebini de EV ile karsilayacagiz ama butceden sikinti yasayabiliriz. Substitutionlar acaba nasil olur vb. Bununla ilgili bazi yorumlar eklemeyi dusunebiliriz.}
}

\subsection{Model Insights} \label{sec:insights}   

\textbf{Service Region and Fleet Sizing Decisions.} We first investigate the profitability of the company under different budget constraints and with and without substitution. The annual net profit, revenue from different trip types, the fixed cost of opening service regions and the cost of relocating cars at the optimal solutions of the models SROP and SROP-S under different budget levels are given in Table \ref{tab:budget}. Note that the demand that can be satisfied by the company depends on the service region opening decisions. In other words, the demand of a potential service region will be lost if a service region is not opened there. Hence, to increase the profit the company should cover more demand by opening more service regions. On the other hand, opening a service region is not sufficient to cover the demand by itself, since there should be sufficiently many vehicles located to that region to satisfy that demand. Therefore, both the total fixed cost of opening service regions and the net profit increase with the budget of fleet sizing, and this can be seen in Table \ref{tab:budget} for both models. When the budget is multiplied by 2 ($B=2M$ vs $B=4M$), the net profit increases approximately by a multiple of 5 in both models. 
\textcolor{black}{Hence, if it is possible, the company should allocate more budget for building the fleet to increase the annual net profit, and the payback period representing the number of years to cover the initial cost of the investment will be shorter for larger initial budget levels. This period reduces further %the numbers are better 
in SROP-S, demonstrating the value of substitution.}

\begin{table}[H]
\centering
\caption{Revenue and Cost Values of SROP and SROP-S Under Different Budget Levels}
\label{tab:budget}
\resizebox{\textwidth}{!}{%
\begin{tabular}{cccccccc}
\hline
                        & Budget & \multicolumn{3}{c}{Revenue(\$K)} & \multicolumn{2}{c}{Cost(\$K)} & Net         \\ \cline{3-7}
Model                   & (\$M)  & Total    & One-Way  & Round-Trip & Fixed        & Relocation     & Profit(\$K) \\ \hline
\multirow{5}{*}{SROP}   & 2      & 1672.06  & 405.15   & 1266.91    & 1447.50      & 52.15          & 172.41      \\
                        & 2.5    & 2200.87  & 648.68   & 1552.19    & 1807.50      & 74.78          & 318.61      \\
                        & 3      & 2744.33  & 949.23   & 1795.10    & 2160.00      & 89.70          & 494.63      \\
                        & 3.5    & 3675.00  & 1627.61  & 2047.39    & 2850.00      & 96.13          & 728.87      \\
                        & 4      & 4253.78  & 1938.98  & 2314.79    & 3210.00      & 109.50         & 934.31      \\ \hline
\multirow{5}{*}{SROP-S} & 2      & 1715.72  & 411.71   & 1304.01    & 1447.50      & 36.82          & 231.40      \\
                        & 2.5    & 2565.25  & 993.59   & 1571.66    & 2115.00      & 39.10          & 411.17      \\
                        & 3      & 3164.92  & 1315.66  & 1849.26    & 2505.00      & 47.22          & 612.70      \\
                        & 3.5    & 3771.12  & 1701.05  & 2070.07    & 2850.00      & 58.78          & 862.34      \\
                        & 4      & 4358.28  & 2018.46  & 2339.83    & 3210.00      & 66.78          & 1081.53     \\ \hline
\end{tabular}%
}
\end{table}

From Table \ref{tab:budget}, we see that the revenue from round-trip arcs is higher than the one-way trips in all settings. But note that the difference decreases with the budget as more vehicles are purchased and more service regions are opened. For instance, in SROP, while the return from round-trips is almost three times larger than the revenue from one-way trips when $B=2M$, they are very close when $B=4M$. Hence, the returns from different trip types depend on the fleet sizing and service region decisions. Moreover, the relocation costs are very small compared to the other cost and return components, but they also increase with the budget due to the same reasons.  

\textbf{Value of Substitution.} We next discuss the effect of substitution by comparing the detailed analysis of the optimal values of the models SROP and SROP-S under different budget $B$ levels. Note that, SROP-S is a generalization of SROP since any feasible solution for the latter is also feasible for the first. Hence, the optimal value of SROP-S should be greater than or equal to the optimal value of SROP, and the difference between the optimal values of these models can be regarded as the value of substitution.
\ignore{
burada fixed ve parking ayri, ama biz modelde hepsine birden fixed diyoruz, o yuzden o iki sutunu birlestirdim.
bir de SROP ve SROP-S icin sonuclarini sirasini degistirdim.
\begin{table}[H]
\centering
\caption{Revenue and Cost Values of SROP and SROP-S Under Different Budget Values}
\label{Table3}
\resizebox{\textwidth}{!}{%
\begin{tabular}{ccccccccc}
\hline
                      &       &           & \multicolumn{3}{c}{Revenue(\$K)} & \multicolumn{3}{c}{Cost(\$K)}  \\ \cline{4-9} 
Budget                & Model & Objective & Total    & One-Way  & Round-Trip & Fixed   & Parking & Relocation \\ \hline
\multirow{2}{*}{3M}   & SROP   & 494.63    & 2744.33  & 949.23   & 1795.10    & 1800.00 & 360.00  & 89.70      \\
                      & SROP-S   & 612.70    & 3164.92  & 1315.66  & 1849.26    & 2100.00 & 405.00  & 47.22      \\
\multirow{2}{*}{3.5M} & SROP   & 728.87    & 3675.00  & 1627.61  & 2047.39    & 2400.00 & 450.00  & 96.13      \\
                      & SROP-S   & 862.34    & 3771.12  & 1701.05  & 2070.07    & 2400.00 & 450.00  & 58.78      \\ \hline
\end{tabular}%
}
\end{table}
}

From Table \ref{tab:budget}, we observe that the total net profit increases with subsitution. For instance, the net profit increases by $34.21\%$ and $29.05\%$ due to substitution when the budget is $B = 2M$ and $B = 2.5M$, respectively. Note that the total service region opening costs are the same for SROP and SROP-S when the budget is $B \in \{2M, 3.5M, 4M\}$, and it is larger for SROP-S in the other budget levels. 

In Table \ref{tab:budget}, we see that the main difference in the revenue is due to the one-way trips covered by these two models especially when the budget is at a medium level. For instance, the revenue obtained from one-way trips increases by 53\% and 39\%  due to substitution when the budget is $B=2.5M$ and $B=3M$, respectively. 
%This increase occurs even though the total fixed cost of SROP-S is larger than that of SROP (6 and 7 regions are opened by the models SROP and SROP-S, respectively). Furthermore, the main difference in the revenue is due to the one-way trips covered by these two models. The revenue obtained from one-way trips increases by $38.60\%$ when substitution is allowed. 

Comparing the results for $B=3M$ and $B=3.5$ in Table \ref{tab:budget} reveals that the value of substitution decreases with the budget. Note that if there exists no budget, carbon emission and capacity constraints, the ideal solution would be to satisfy all demand by its own car type. Hence, when the budget for fleet sizing is larger, since more cars will be available in the service regions, satisfying the demand by the actual car type demanded will occur more, and the value of substitution will decrease. But, as it can be observed from Table \ref{tab:budget}, though the effect of substitution is smaller compared to $B=3M$, the total net profit is increased by $ 18.31\%$ due to substitution when the budget is $B=3.5M$. Note that the total service region opening costs are the same for the models in these budget levels.

We also observe an interesting result from Table \ref{tab:budget} where the total relocation cost of SROP-S is lower than that of SROP under all budget levels. We note that substitution between different car types works as a relocation in SROP-S. In other words, instead of using a worker to relocate a vehicle between two regions (and observing a cost), using that vehicle for satisfying the demand for the other vehicle type (and getting a return) helps the company to balance the cars at the service regions, and this reduces the relocation costs and increases the return. 

\begin{table}[H]
\centering
\caption{Value of Substitution and Average Flows for $B=\$3M$} \label{tab:flowsbudget3}
\resizebox{\textwidth}{!}{%
\begin{tabular}{ccccccccc}
\hline
                   &           & \multicolumn{4}{c}{Average Flow}          & Subs.             & Objective         & Achieved           \\ \cline{3-6}
Penalty ($p$)           & Commodity & One-way & Round-Trip & Relocation & Idle  & Rate(\%) &  Increase(\%)             & Increase(\%)           \\ \hline
\multirow{2}{*}{$\rightarrow \infty$ }    & E-E & 66.74 & 70.96 & 15.55 & 80.90 & \multirow{4}{*}{-} & \multirow{4}{*}{-}     & \multirow{4}{*}{-}     \\
                   & E-G       & -       & -          & -          & -     &          &                        &                        \\
 \multirow{2}{*}{($E:48, G: 48$)}    & G-G       & 67.12   & 70.78      & 15.17      & 82.46 &          &                        &                        \\
                   & G-E       & -       & -          & -          & -     &          &                        &                        \\ \hline
\multirow{2}{*}{4} & E-E       & 81.52   & 63.28      & 7.61       & 48.94 & -        & \multirow{4}{*}{19.00} & \multirow{4}{*}{64.85} \\
                   & E-G       & 7.97    & 4.12       & -          & -     & 5.57     &                        &                        \\
 \multirow{2}{*}{($E:46, G: 53$)}    & G-G       & 86.21   & 69.07      & 8.97       & 59.99 & -        &                        &                        \\
                   & G-E       & 10.64   & 7.91       & -          & -     & 8.60     &                        &                        \\ \hline
\multirow{2}{*}{2} & E-E       & 80.25   & 62.90      & 7.46       & 47.97 & -         & \multirow{4}{*}{23.88} & \multirow{4}{*}{81.50} \\
                   & E-G       & 9.12    & 5.18       & -          & -     & 6.59     &                        &                        \\
 \multirow{2}{*}{($E:46, G: 53$)}     & G-G       & 84.55   & 68.17      & 8.71       & 59.07 & -        &                        &                        \\
                   & G-E       & 11.77   & 8.97       & -          & -     & 9.61     &                        &                        \\ \hline
\multirow{2}{*}{$\rightarrow 0$} & E-E & 79.77 & 61.97 & 7.18  & 48.83 & -                   & \multirow{4}{*}{29.30} & \multirow{4}{*}{100}\\
                   & E-G       & 10.02   & 5.86       & -          & -     & 7.31     &                        &                        \\
 \multirow{2}{*}{($E:46, G: 53$)}   & G-G       & 84.59   & 66.94      & 8.90       & 59.66 & -        &                        &                        \\
                   & G-E       & 13.36   & 9.68       & -          & -     & 10.69    &                        &                        \\ \hline
\end{tabular}%
}
\end{table}

In Table \ref{tab:flowsbudget3}, we present the average number of one-way and round-trips, relocation actions, and idle vehicles waiting in the service regions under different substitution penalty prices $p$ when the budget is $B=3M$. Note that $p \rightarrow \infty$ represents SROP since SROP-S reduces to SROP when $p\rightarrow \infty$ as substitution will not be used in this case. For the other extreme case, where substitution is allowed and not penalized, $p\rightarrow 0$, %we do not consider the setting $p=0$ since in that case there exists no distinction between the car types. Instead, 
we consider a very small but positive $p$ value ($p=0.001$) to observe meaningful results for different commodities. Considering two car types, we have four commodities in SROP-S, and the commodity $a-b$ represents the case where the demand for car type $b$ is satisfied by car type $a$, for $a, b \in \{E,G\}$. Note that the company gains money over the flows on one-way and round-trip arcs, loses money due to the flows on relocation arcs, and has no gain or cost on the flows on idle arcs (though they might be also perceived as loss on potential profit).  

From Table \ref{tab:flowsbudget3}, we observe that the average number of idle vehicles and the relocation actions decrease dramatically with substitution. Note that the average number of idle vehicles is around 80 for both car types in SROP while it is less than 50 and 60 for electric and gasoline cars, respectively, in SROP-S. Similarly, the average number of relocation movements are around 15 and 8 in SROP and SROP-S, respectively. Hence, substitution provides the company a flexibility for eliminating the non-value adding operations (relocation and being idle in our case) and increasing the net profit. 
Moreover, as it can be seen in Table \ref{tab:flowsbudget3}, one-way trips are preferred more  when substitution is allowed since one-way trips also serve as relocation and their unit profit is larger. Additionally, average number of one-way trips increases while the average number of round-trip trips decreases in SROP-S compared to SROP. 

When we compare the results of SROP-S for different penalty $p \in \{0,2,4\}$ values, we observe that though the average flows of substitution commodities ($E-G$ and $G-E$) decrease, the average number of trips do not change so much for each commodity. Hence, we can say that the solution of SROP-S is robust with respect to the different values of $p$ except the case $p \rightarrow \infty $ which represents SROP. 

In the last three columns of Table \ref{tab:flowsbudget3}, we present the percentage of substitution, increase in the net profit and the relative increase in the net profit for different $p$ values. As it is expected, the substitution rates and the net profit increases when $p$ is decreasing. Moreover, the substitution rates are larger for gasoline vehicles under all settings. Since gasoline vehicle is cheaper, when substitution is allowed, the company buys more gasoline cars (the numbers of electric and gasoline cars purchased are given under the column Penalty in parenthesis) and uses them for satisfying the demand for electric cars. To calculate Achieved Increase rates for each penalty parameter, we take the ratio of the percentage increase in the net profit under that penalty parameter to the increase in the net profit when $p \rightarrow 0$. For instance, though the net profit is increased by 19\% when $p=4$ compared to $p \rightarrow \infty$, this increase actually corresponds to the 64.85\% of the maximum possible increase in the net profit, showing that a significant amount of increase is achieved with the substitutions.   %Hence, the solution and the net profit do not change so much for small values of $p$. 
Appendix \ref{sec:AppendixValueofSubstitution} provides additional results on the value of substitution when the car sharing company has higher budget. 

\begin{figure}[H]
  \caption{Visual Comparison of Demand Satisfaction Rates for $B=3M$ (left) and for $B=3.5M$ (right)}
 \begin{minipage}{0.5\linewidth}
  \centering
  \includegraphics[scale=0.6]{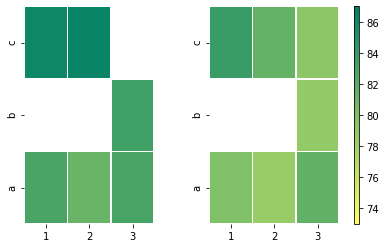}
  
  SROP (left) and SROP-S (right)
  \label{fig:visual3}
 \end{minipage}
 \begin{minipage}{0.5\linewidth}
  \centering  \includegraphics[scale=0.6]{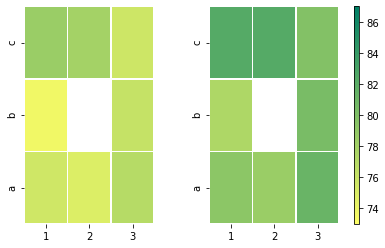}
  
    SROP (left) and SROP-S (right) \label{fig:visual3.5}
 \end{minipage}
\end{figure}

\textbf{Demand Satisfaction Levels.} We next discuss the demand satisfaction rates for these two models under two different budget levels. Notice that the demand that can be covered strongly depends on the service region opening decisions. Hence, in Table \ref{tab:demandsat} we report the demand satisfaction rates with respect to two different values. We first present the number of service regions opened under the column SR. In columns 4 and 5, we give the percentage of the satisfied demand for electric (EV) and gasoline vehicles with respect to the total demand of the service regions opened. In columns 6 and 7, we report these percentages with respect to the total demand of the whole system. In the remaining columns, the percentage of the satisfied demand for each region is given. Positions of these regions, i.e. 1-a, 2-a, etc., with respect to each other can be seen in Figure \ref{fig:visual3}, %and \ref{fig:visual3.5} 
where the demand satisfaction rates for the regions are illustrated visually with different colors.  

\begin{table}[htbp]
\centering
\caption{Demand Satisfaction Rates for SROP and SROP-S Under Different Budget Parameters}\label{tab:demandsat}
\resizebox{\textwidth}{!}{%
\begin{tabular}{cccccccccccccccc}
\hline
       &       &    & \multicolumn{2}{c}{SR(\%)} & \multicolumn{2}{c}{Total(\%)} & \multicolumn{9}{c}{Service Region Based Demand Satisfaction(\%)}                                             \\ \cline{4-16} 
Budget & Model & SR & EV        & Gasoline    & EV       & Gasoline    & 1-a     & 2-a     & 3-a     & 1-b     & 2-b    & 3-b     & 1-c     & 2-c     & 3-c     \\ \hline
\multirow{2}{*}{3M}   & SROP & 6 & 81.57 & 81.37 & 42.53 & 42.25 & 82.94 & 81.24 & 83.02 & 0.00  & 0.00 & 83.43 & 86.23 & 86.5  & 0.00  \\
       & SROP-S   & 7  & 76.28     & 77.25      & 50.57     & 51.13      & 79.98 & 78.61 & 81.59 & 0.00  & 0.00 & 78.93 & 83.85 & 81.48 & 79.26 \\ \cline{2-16}
\multirow{2}{*}{3.5M} & SROP & 8 & 72.19 & 74.24 & 59.35 & 60.76 & 75.63 & 74.96 & 77.03 & 73.81 & 0.00 & 76.20 & 78.53 & 77.98 & 75.75 \\
       & SROP-S   & 8  & 76.01     & 77.08      & 62.50     & 63.10      & 79.32 & 78.57 & 81.24 & 77.48 & 0.00 & 80.37 & 82.3  & 82.35 & 79.80 \\ \hline
\end{tabular}%
}
\end{table}
\vspace{-8mm}
\vspace{-5mm}
%\begin{figure}[h] %htbp
%    \centering
%    \caption{Visual Comparison of Demand Satisfaction Rates of SROP (left) and SROP-S (right) for $B=3M$} \label{fig:visual3}
%    \includegraphics[scale=0.6]{heatmap1a.png}\hfill
    %width=.8\textwidth,height=.25\textheight]
%\end{figure}
From Table \ref{tab:demandsat}, we see that SROP has higher demand satisfaction levels for both car types compared to SROP-S when $B=3M$ and only the demand of opened service regions are considered. On the other hand, the demand satisfaction rates are higher for SROP-S with respect to the total demand of the system under the same budget level. Note that the number of service regions opened are not the same for SROP and SROP-S when $B=3M$. Hence, when $B=3M$, SROP-S opens one more service region and this reduces the demand satisfaction rates of the other service regions opened (see the last 9 columns). This shows the trade-off between opening a new service region and satisfying more demand by serving less regions under a limited budget for fleet sizing. SROP-S allocates some of the purchased cars to the additional region it opens, region 3-c, and due to this fact, the demand satisfaction rates in the other opened regions decrease compared to SROP. But note that this decision increases the total net profit of the company by $23.87\%$ (see Table \ref{tab:budget}). When $B=3.5M$, since the number of opened regions is the same for both models, the demand satisfaction levels are larger for SROP-S in all opened regions. Figure \ref{fig:visual3} %and \ref{fig:visual3.5} 
illustrates these observations. In Figure \ref{fig:visual3} when $B=3M$, since SROP opens less service regions, demand satisfaction levels are higher (darker) in all \textit{opened} regions. In Figure \ref{fig:visual3.5} when $B=3.5M$, since both models open the same regions, the demand satisfaction rates are larger (darker) in SROP-S in all regions. 
%We have two additional observations from Figures \ref{fig:visual3} and \ref{fig:visual3.5}. First, 
Additionally, note that the regions that are decided to be served by SROP are also covered by SROP-S. Hence, the service region opening decisions are not affected so much by the substitution decision though additional ones might be opened in SROP-S due to the flexibility it provides to the company.

\textbf{Carbon emission.} We report the results of SROP and SROP-S under two different budget levels $B \in \{3M, 3.5M\}$ and two different carbon emission allowances $H \in \{0.3, 0.5\}$ in Table \ref{tab:carbon}. As the problems become more restricted for small $H$ values, the net profit increases with $H$. Note that the value of substitution becomes more prominent when the problem is more restricted, i.e. $H=0.3$ and $B=3M$ with more than 58\% increase in net profit.   

The unitary carbon allowance $H$ directly affects the percentage of the car types in the fleet and this can be observed from Table \ref{tab:carbon}. Although the fleet sizes are very close under different $H$ levels, the ratio of the number of electric and gasoline cars changes significantly. As the gasoline cars are cheaper but have larger emissions, the percentage of gasoline cars in the fleet increases with $H$. This also affects the number of service regions opened. 
Appendix \ref{sec:AppendixLowerCarbon} provides visual comparisons of demand satisfaction rates of SROP and SROP-S models under different budget levels. 
Notice that more service regions are opened by both models when $H=0.3$ compared to $H=0.5$. For smaller $H$ values, since more EVs are purchased and there are separate parking space capacities for electric and gasoline cars, the models prefer to open more service regions to use these electric cars for satisfying more demand. 
%\caption{Planning Decision Effect of Carbon Emission Parameter}
\begin{table}[H]
\centering
\caption{The Effect of Carbon Emission Allowance on Planning Decisions} \label{tab:carbon}
%\footnotesize{}
%\resizebox{\textwidth}{!}{%
\begin{tabular}{cclccccc}
\hline
$B$ & $H$ & Model & Obj.(\$K) & Obj. Inc.(\%) &  SR & \# of Cars (G/E) & Total \# of Cars  \\ \hline
\multirow{4}{*}{3M}  &  \multirow{2}{*}{0.3}   & SROP   & 269.5 & - & 8 & 38/58 & 96         \\
                       &                       & SROP-S & 427.3 & 58.52 & 8 & 38/58 & 96   \\ \cline{2-8}
&\multirow{2}{*}{0.5}   & SROP   & 494.6 & - & 6 & 48/48 & 96         \\
                       &                       & SROP-S & 612.7 & 23.88& 7 & 53/46 & 99   \\   \hline
\multirow{4}{*}{3.5M}  & \multirow{2}{*}{0.3} & SROP   & 461.5 & -& 9 & 44/68 & 112        \\
                     &                       & SROP-S & 645.0 & 39.76& 9 & 44/68 & 112  \\ \cline{2-8}
& \multirow{2}{*}{0.5} & SROP   & 728.9 & -& 8 & 59/56 & 115        \\
                     &                       & SROP-S & 862.4 & 18.32& 8 & 59/56 & 115  \\  \hline
\end{tabular}%
%}
\end{table}
\subsection{Computational Performances} \label{sec:performance}
In the second part of our computational study, we evaluate the performance of our decomposition based solution algorithm for solving large problem instances.   
We consider problem instances with $|W| \in \{50, 100, 200\}$ scenarios under two different budget levels $B \in \{3M, 3.5M \}$, and generate three random problem instances for each setting.  We set the time limit to four hours. %14,440 seconds. 
For both of the models SROP and SROP-S, we compare the Deterministic Equivalent Formulations (DEF) of the models SROP and SROP-S given in Section \ref{sec:ProblemForm}, and the Branch-and-Cut algorithms without (B\&C) and with the enhancements (B\&C+) presented in Section \ref{sec:solution}. In our preliminary experiments, we test all enhancements stated in Section \ref{sec:imp&enhanc}, and observe that the contribution of enhancement (iii) is limited compared to the first two. Hence, we omit enhancement (iii), and use (i) and (ii) in our B\&C+ algorithm.  
All tests are performed on a personal computer running Microsoft Windows 10 64 bit operating system at Intel i5-10210U 1.60 GHz processor with 16 GB RAM at 8 threads. All formulations and algorithms are implemented and solved in Gurobi 9.5.2 and Phyton interface with the default settings. 

We present the results in Table \ref{tab:comp} where ST and Gap represent the solution time (in seconds) and the optimality gap, respectively. If the solver terminates due to the time limit, it is stated as TL under the column ST, and the optimality gap reported by the solver at the end of the time limit  is presented under the column Gap. If no feasible solution could be found by the solver within the time limit, we note it as $N/A$ under the column Gap. For each instance, we emphasize the best result in terms of the solution time and the optimality gap in bold.  

\begin{table}[h]
\centering
\caption{Computational Comparison of Model Types for Different Parameters and Problem Sizes } \label{tab:comp}
\resizebox{\textwidth}{!}{%
\begin{tabular}{ccc|cccccc|cccccc}
\hline
                     &                       &       & \multicolumn{6}{|c|}{SROP}                          & \multicolumn{6}{c}{SROP-S}                         \\ \cline{4-15} 
$|W|$ & $B$ & Ins & \multicolumn{2}{|c}{DEF} & \multicolumn{2}{c}{B\&C} & \multicolumn{2}{c|}{B\&C+} & \multicolumn{2}{c}{DEF} & \multicolumn{2}{c}{B\&C} & \multicolumn{2}{c}{B\&C+} \\
                 & (\$)            & \# & ST    & Gap(\%) & ST   & Gap(\%) & ST   & Gap(\%) & ST    & Gap(\%) & ST   & Gap(\%) & ST    & Gap(\%) \\ \hline
\multirow{6}{*}{50}  & \multirow{3}{*}{3M}   & 1     & 1085  & 0.00    & 1326 & 0.00    & \textbf{1282} & 0.00    & 5681  & 0.00    & 3138 & 0.00    & \textbf{2961}  & 0.00    \\
                     &                       & 2     & 1513  & 0.00    & \textbf{1478} & 0.00    & 1537 & 0.00    & 6208  & 0.00    & 3768 & 0.00    & \textbf{2457}  & 0.00    \\
                     &                       & 3     & \textbf{719}   & 0.00    & 1359 & 0.00    & 1258 & 0.00    & 5314  & 0.00    & 3974 & 0.00    & \textbf{3483}  & 0.00    \\ \cline{2-15} 
                     & \multirow{3}{*}{3.5M} & 1     & 690   & 0.00    & 653  & 0.00    & \textbf{510}  & 0.00    & 2731  & 0.00    & 1442 & 0.00    & \textbf{593}   & 0.00    \\
                     &                       & 2     & 692   & 0.00    & 551  & 0.00    & \textbf{317}  & 0.00    & 3191  & 0.00    & 1261 & 0.00    & \textbf{770}   & 0.00    \\
                     &                       & 3     & 640   & 0.00    & 485  & 0.00    & \textbf{427}  & 0.00    & 3525  & 0.00    & \textbf{1291} & 0.00    & 1325  & 0.00    \\ \cline{1-15} 
\multirow{6}{*}{100} & \multirow{3}{*}{3M}   & 1     & 4093  & 0.00    & 2555 & 0.00    & \textbf{2156} & 0.00    & TL    & 100.00  & 6379 & 0.00    & \textbf{6076}  & 0.00    \\
                     &                       & 2     & 3823  & 0.00    & \textbf{2817} & 0.00    & 2886 & 0.00    & TL    & 100.00  & 6451 & 0.00    & \textbf{6167}  & 0.00    \\
                     &                       & 3     & 3080  & 0.00    & 2449 & 0.00    & \textbf{2281} & 0.00    & TL    & 39.60   & 7053 & 0.00    & \textbf{6532}  & 0.00    \\ \cline{2-15} 
                     & \multirow{3}{*}{3.5M} & 1     & 2486  & 0.00    & 1144 & 0.00    & \textbf{862}  & 0.00    & 10997 & 0.01    & 2403 & 0.00    & \textbf{1310}  & 0.00    \\
                     &                       & 2     & 2368  & 0.00    & \textbf{1318} & 0.00    & 1336 & 0.00    & 12203 & 0.01    & 2974 & 0.00    & \textbf{2440}  & 0.00    \\
                     &                       & 3     & 2417  & 0.00    & 1039 & 0.00    & \textbf{988}  & 0.00    & TL    & 0.00    & 3298 & 0.00    & \textbf{2558}  & 0.00    \\ \hline
\multirow{6}{*}{200} & \multirow{3}{*}{3M}   & 1     & TL    & 35.31   & 4444 & 0.00    & \textbf{3507} & 0.00    & TL    & N/A     & TL   & 0.02    & \textbf{11175} & 0.00    \\
                     &                       & 2     & TL    & 100.00  & 6859 & 0.00    & \textbf{5236} & 0.00    & TL    & 100.00  & TL   & 3.26    & \textbf{10120} & 0.00    \\
                     &                       & 3     & 14110 & 0.00    & \textbf{5085} & 0.00    & 6185 & 0.00    & TL    & 100.00  & TL   & 17.43   & TL &\textbf{ 0.69 }   \\ \cline{2-15} 
                     & \multirow{3}{*}{3.5M} & 1     & 11192 & 0.00    & 3115 & 0.00    & \textbf{1403} & 0.00    & TL    & 73.89   & 5476 & 0.00    & \textbf{3810}  & 0.00    \\
                     &                       & 2     & 12404 & 0.00    & 2171 & 0.00    & \textbf{1470} & 0.00    & TL    & N/A     & 4963 & 0.00    & \textbf{2980}  & 0.00    \\
                     &                       & 3     & 10332 & 0.00    & 1838 & 0.00    & \textbf{1628} & 0.00    & TL    & 80.06   & \textbf{6140} & 0.00    & 6317  & 0.00    \\ \hline
\end{tabular}%
}
\end{table}

As it can be seen from Table \ref{tab:comp}, the single stage formulations DEF can be solved to optimality within the time limit only for small problem instances. The gaps reported for DEF demonstrate that it is hard to obtain a good quality solution for large problem instances using DEF. On the other hand, the decomposition based algorithms B\&C and B\&C+ perform better than DEF almost in all problem instances. Both algorithms solve SROP instances within the time limit, but the solution times are mostly better for B\&C+. 
The size of SROP-S is larger than SROP for a given instance since the number of commodities in these models are four and two, respectively. Therefore, it is harder to solve SROP-S, in general, and this can be observed from Table \ref{tab:comp}. Again, the best algorithm for solving SROP-S is B\&C+, which means that the  enhancements presented in Section \ref{sec:solution} improves the performance of the algorithm B\&C. Moreover, the problems become relatively easier to solve for larger $B$ levels.

\section{Conclusion}\label{sec:Conclusion}
In this paper, we study the service region design and operational planning problem for a car sharing company that constructs a mixed fleet of gasoline and EVs under the demand uncertainty for one-way and round-trip rental requests and the carbon emission constraints. We introduce substitution to the problem by allowing satisfying the demand for a specific car type using the other type. We develop a two-stage stochastic mixed-integer programming model for the problem, and propose an exact decomposition based algorithm. Our computational experiments reveal the success of our algorithm in solving larger problem instances compared to the off-the-shelf solver. We discuss the impact of different problem parameters on the  decisions of the problem through a case study based on real data sets. The results of the case study indicate the value of substitution for both increasing the profit of the company and the demand satisfaction level of the customers. Since substitution gives a flexibility to the car sharing company, we observe that less rebalancing operations are required when substitution is allowed. Besides, the results show the important effect of the fleet sizing budget and the carbon emission allowance on the service region opening and the fleet allocation decisions. 

As future research directions, the model considered in this paper can be extended from several directions. %has several limitations. 
First, we assume that the customers accept substitution independent from the car type that is substituted and the price offered. In a future research, substitution should be studied in more detail by including the customer behaviour and also the pricing strategy of the company for substitution. Second, in this study we approximate the demand for the next $T$ periods using two-stage stochastic programming, and a more appropriate approach might be to formulate the problem as a multi-stage stochastic program which is more challenging to solve. 

%Daha fazla sey soylenebilir mi?

\ignore{
\section{Notes}

\begin{itemize}
    \item Demand satisfaction tablosundaki degerler icin normalize hallerini yeni bir tablo yaparak appendixe ekleyebiliriz
    \item 2. modelin aciklamasini yaparken commodity aciklamalarini vermeye gerek yok, genel olarak verelim; sonra su sekilde olabilir diye bir ornek model eklenebilir
    \item ikinci modelin dualinde relocation kisidini sormayi unutma!
\end{itemize}

}

%\section{Appendix}

\bibliographystyle{informs2014}
\bibliography{references} % if more than one, comma separated

\begin{APPENDICES}

\section{Proofs}

\subsection{Proof of Proposition \ref{prop:TU-SROP}}
\label{proof:TU-SROP}

We observe that constraints \eqref{eq:M1_11} represent the flow conservation constraints, corresponding to a node-arc incidence matrix, which is totally unimodular. Since the capacity constraints \eqref{eq:M1_12} - \eqref{eq:M1_15} define a unit row for each round-trip arc, two unit rows for each one-way arc, and a unit row for each idle arc, we can ignore them. Then, the result follows. 
%Flow variable $y_{akw}$ stated as an integer variable, however, since the second stage problem is a network flow problem with flow conversation and capacity constraints, it can be written as a totally unimodular matrix. Thus, we can consider the flow variable as a continuous variable as well.

\subsection{Proof of Proposition \ref{prop:TU-SROP-S}}
\label{proof:TU-SROP-S}

Consider the constraint matrix of \eqref{eq:M2_1} -\eqref{eq:M2_5} for a given scenario $w \in W$. 
%Since the variables $y_{alw}$ are defined for each arc and commodity pair $(a,l)$, we can assume that the nodes and arcs are copied for car type and commodity type, respectively. Hence, 
Observe that the balance constraints \eqref{eq:M2_1} define a node-arc incidence matrix. The remaining constraints \eqref{eq:M2_2}-\eqref{eq:M2_5} partition the arc set into three groups: \eqref{eq:M2_2} is written only for round-trip arcs, \eqref{eq:M2_3} and \eqref{eq:M2_4} are for one-way arcs and \eqref{eq:M2_5} is for idle arcs. Since \eqref{eq:M2_5} defines a unit row for each idle arc, we can ignore them. Moreover, both $F_k$ and $\hat{F}_k$ partition the commodity set $L$ into $k$ subsets. In other words, a variable $y_{alw}$ is seen only in one of the constraints \eqref{eq:M2_2} if $a \in A^{two}$. Similarly, $y_{alw}$ is seen in one of the constraints \eqref{eq:M2_3} and in one of the \eqref{eq:M2_4} if $a$ is $a \in A^{one}$. 

Let $G_1, G_2, G_3, G_4$ be the set of rows due to constraints \eqref{eq:M2_1}, \eqref{eq:M2_2}, \eqref{eq:M2_3}, and \eqref{eq:M2_4}, respectively. We will prove that for any subset $R$ of rows of the constraint matrix, there exists a partition $R_1$, $R_2$ of $R$ such that the difference between these subsets is $\{0,1,-1\}$ in any column.  For any $R$, we construct our rule for partition as follows:
\begin{itemize}
    \item All rows in $R \cap G_1$ will be included in $R_1$. %Then the summation of the rows in $R_1$ will be $\{0,-1,1\}$ for each column.
    \item Each row in $R \cap G_2$ corresponds to an arc $a=(n_{it},n_{is}) \in A^{two}$ and car type $k \in K$: 
    \begin{itemize}
        \item If $R \cap G_1$ includes the rows for $i,t,k$ and $i,s,k$, then we can include the row for $a,k$ to any of $R_1$ and $R_2$. 
        \item If $R \cap G_1$ includes the row for $i,t,k$ but not for $i,s,k$, then we add the row for $a,k$ to $R_1$. 
        \item If $R \cap G_1$ includes the row for $i,s,k$ but not for $i,t,k$, then we add the row for $a,k$ to $R_2$. 
    \end{itemize}
    \item Each row in $R \cap G_3$ and $R \cap G_4$ corresponds to an arc $a=(n_{it},n_{js}) \in A^{one}$ and car type $k \in K$. 
    \begin{itemize}
        \item If $R \cap (G_3 \cup G_4)$ includes two rows for $a,k$, then we can include one of the rows to $R_1$ and the other one to $R_2$. 
        \item If $R \cap (G_3 \cup G_4)$ includes only one row for $a,k$, then
        \begin{itemize}
            \item if  $R \cap G_1$ includes both of the rows for $i,t,k$ and $j,s,k$, then we can include the row for $a,k$ to any of $R_1$ and $R_2$.
            \item if $R \cap G_1$ includes the row for $i,t,k$ but not for $j,s,k$, then we add the row for $a,k$ to $R_1$. 
            \item if $R \cap G_1$ includes the row for $j,s,k$ but not for $i,t,k$, then we add the row for $a,k$ to $R_2$. 
        \end{itemize}
    \end{itemize}
\end{itemize}
Note that the partition of rows in $R \cap G_2$ and $R \cap (G_3 \cup G_4)$ are independent from each other since they include different columns, i.e. $G_2$ includes the variables defined for two-way arcs while $G_3 \cup G_4$ includes the variables for one-way arcs. Hence, due to the construction of $R_1$ and $R_2$, the difference between the summation of rows in $R_1$ and $R_2$ will be in $\{0,-1,1\}$ for any column (variable), and the result follows. 

%Simdilik tez formati icin buraya almadim, ancak 32 sayfa siniri sebebiyle sonra prooflari buraya aliriz.  
\section{Additional Models and Instance Information}

\subsection{Dual of the Second Stage Problem for SROP-S}
\label{sec:DualSecondStageSROPS}
Given the optimal solution $(\hat{\mathbf{z}},\hat{\mathbf{x}}, \hat{\mathbf{\bar{Q}}})$ for RMP, the dual of the second-stage problem \eqref{eq:SROP-S-TwoStage-2nd} for each scenario $w \in W$ can be written in \eqref{eq:SROP-S-Benders-dualsub}. 
\begin{subequations}
\label{eq:SROP-S-Benders-dualsub}
\begin{alignat}{1}
\begin{split}
\displaystyle \min \quad &  \sum_{i\in I} \sum_{k\in K} \hat{x}_{ik} (\beta_{i0kw} - \beta_{iTkw})  + \sum_{i \in I}\sum_{a \in  A^{two}(i)}\sum_{k \in K} u_{akw} \hat{z}_{i} \alpha_{akw}  + 
 \sum_{i \in I}\sum_{a \in  A^{one}(i^{+})} 
 \sum_{k \in K} u_{akw} \hat{z}_{i} \gamma^{1}_{akw} \\ 
& +  \sum_{i \in I}\sum_{a \in  A^{one}(i^{-})}\sum_{k \in K} u_{akw} \hat{z}_{i} \gamma^{2}_{akw} +  
     \sum_{i \in I} \sum_{k \in K} C_i^k \hat{z}_{i} (\sum_{a \in A^{idle}(i)} \lambda_{akw})
\end{split}\\
\text{s.t.} \quad
& \beta_{itk^{1}w} - \beta_{jsk^{1}w} + \gamma^{1}_{ak^2w}+ \gamma^{2}_{ak^2w}  \geq r_{al} \quad  \forall a = (n_{it},n_{js}) \in A^{one}, \> \forall k^1,k^2 \in K, \; l \in F^{k^1} \cap \hat{F}^{k^2}, \\
& \beta_{itk^{1}w} - \beta_{isk^{1}w} +  \alpha_{ak^2w}  \geq r_{al} \quad \forall a = (n_{it},n_{is}) \in A^{two},  \> \forall k^1,k^2 \in K, \; \> l \in F^{k^1} \cap \hat{F}^{k^2}, \\
& \beta_{itkw} - \beta_{j,t+\zeta_{ij},kw} \geq r_{al} \quad \forall a = (n_{it},n_{j,t+\zeta_{ij}}) \in A^{rel}, \> \forall k \in K, \> l \in \in F^{k} \cap \hat{F}^{k}, \\
& \beta_{itkw} - \beta_{i,t+1,kw} + \lambda_{akw} \geq 0 \quad \forall a = (n_{it},n_{i,t+1}) \in A^{idle}, \> \forall k \in K,   \\
&\gamma^{1}_{akw}, \gamma^{2}_{akw} \geq 0 \quad \forall a \in  A^{one}, \> \forall k \in K, \\
& \alpha_{akw} \geq 0 \quad \forall a \in  A^{two}, \> \forall  k \in K, \\
& \lambda_{akw} \geq 0 \quad \forall a \in A^{idle}, \> \forall k \in K.
\end{alignat}
\end{subequations}

\subsection{Service Region Capacities and Costs}
\label{sec:ServiceRegionCapacities}
\begin{table}[H]
\centering
\caption{Capacities and fixed costs for each service region} 
\label{tab:data1}
\begin{tabular}{c|ccccccccc} \hline
Region ($i$) &  1 & 2 & 3 & 4 & 5 & 6 & 7 & 8 & 9 \\ \hline
Capacity for each car type($C_i^k$) & 6 & 9 & 7 & 6 & 8 & 9 & 8 & 9 & 6\\ 
Fixed cost ($f_i$)($\$K$) & 345 & 367.5 & 352.5 & 345 & 360 & 367.5 & 360 & 367.5 &  345 \\ \hline
 \end{tabular}
 \end{table}

\section{Additional Results}
% \section{<Title of Section B>}
% etc
\subsection{Value of Substitution under Higher Budget}
\label{sec:AppendixValueofSubstitution}

Table \ref{tab:flowsbudget3.5} provides the value of substitution when budget $B=\$3.5M$. Compared to Table \ref{tab:flowsbudget3}, number of cars purchased increases due to the increase in budget from $B=\$3M$ to $B=\$3.5M$. Thus, there is less need for substitution, resulting in slightly less  increases in net profit values. Similar to the previous results, majority of maximum potential return is achieved under various penalty parameter values. Additionally, car flows are robust to the changes in penalty parameters. 

\begin{table}[H]
\label{tab:flowsbudget3.5}
\centering
\caption{Value of Substitution and Average Flows for $B=\$3.5M$}
\resizebox{\textwidth}{!}{%
\begin{tabular}{cccccccccc}
\cline{1-9}
                                      &           & \multicolumn{4}{c}{Average Flow}          & Subs.              & Objective              & Achieved               & \multicolumn{1}{l}{} \\ \cline{3-6}
Penalty                               & Commodity & One-way & Round-Trip & Relocation & Idle  & Rate(\%)           & Increase(\%)           & Return(\%)             & \multicolumn{1}{l}{} \\ \cline{1-9}
\multirow{2}{*}{$\rightarrow \infty$} & E-E       & 114.69  & 77.53      & 15.26      & 61.70 & \multirow{4}{*}{-} & \multirow{4}{*}{-}     & \multirow{4}{*}{-}     & \multicolumn{1}{l}{} \\
                                      & E-G       & -       & -          & -          & -     &                    &                        &                        &                      \\
\multirow{2}{*}{($E:56, G: 59$)}      & G-G       & 117.39  & 81.04      & 17.66      & 69.28 &                    &                        &                        &                      \\
                                      & G-E       & -       & -          & -          & -     &                    &                        &                        &                      \\ \cline{1-9}
\multirow{2}{*}{4}                    & E-E       & 110.14  & 73.14      & 9.52       & 55.35 & -                  & \multirow{4}{*}{14.46} & \multirow{4}{*}{63.89} & \multicolumn{1}{l}{} \\
                                      & E-G       & 10.31   & 5.98       & -          & -     & 6.07               &                        &                        &                      \\
\multirow{2}{*}{($E:56, G: 59$)}      & G-G       & 114.09  & 75.68      & 10.90      & 59.59 & -                  &                        &                        &                      \\
                                      & G-E       & 12.51   & 7.06       & -          & -     & 7.17               &                        &                        &                      \\ \cline{1-9}
\multirow{2}{*}{2}                    & E-E       & 108.24  & 72.35      & 9.50       & 54.02 & -                  & \multirow{4}{*}{18.32} & \multirow{4}{*}{80.93} & \multicolumn{1}{l}{} \\
                                      & E-R       & 12.05   & 6.92       & -          & -     & 7.08               &                        &                        &                      \\
\multirow{2}{*}{($E:56, G: 59$)}      & G-G       & 112.21  & 74.98      & 10.63      & 59.22 & -                  &                        &                        &                      \\
                                      & G-E       & 13.82   & 8.17       & -          & -     & 8.21               &                        &                        &                      \\ \cline{1-9}
\multirow{2}{*}{$\rightarrow 0$}                & E-E       & 108.01  & 71.08      & 9.50       & 55.38 & -                  & \multirow{4}{*}{22.62} & \multirow{4}{*}{100} & \multicolumn{1}{l}{} \\
                                      & E-G       & 13.35   & 7.70       & -          & -     & 7.85               &                        &                        &                      \\
\multirow{2}{*}{($E:56, G: 59$)}      & G-G       & 111.77  & 73.73      & 10.54      & 59.12 & -                  &                        &                        &                      \\
                                      & G-E       & 15.46   & 9.23       & -          & -     & 9.23               &                        &                        &                      \\ \cline{1-9}
\end{tabular}%
}
\end{table}

\subsection{Demand Satisfaction Levels under Lower Carbon Parameter Value}
\label{sec:AppendixLowerCarbon}

Figure \ref{fig:visual3_LowerCarbon} %and \ref{fig:visual3.5_LowerCarbon} 
provides which regions are open and how demand is satisfied under different budget levels when $H = 0.3$. Since the problem becomes more restrictive due to lower carbon parameter value of $H = 0.3$, demand satisfaction rates are smaller compared to the cases when $H = 0.5$. Different than Figure \ref{fig:visual3} under $B=3M$, %and \ref{fig:visual3.5}, 
which has different open service regions under SROP and SROP-S models, when the model becomes more restrictive in terms of carbon allowance, then both models open the same set of service regions under both budget levels in Figure \ref{fig:visual3_LowerCarbon}. %and \ref{fig:visual3.5_LowerCarbon}. 
On the other hand, demand satisfaction percentages are higher under the SROP-S model in all open service regions by utilizing substitution in satisfying customer demand. 

%\begin{figure}[h]
%    \centering
%    \caption{Visual Comparision of Demand Satisfaction for H=0.3 and B=3M}
%    \label{fig:visual3_LowerCarbon}
    %\includegraphics[width=.8\textwidth,height=.25\textheight]{H03-heatmap1-v2.png}\hfill
%\end{figure}
%\begin{figure}[h]
%    \centering
%    \caption{Visual Comparision of Demand Satisfaction for H=0.3 and B=3.5M}
    %\label{fig:visual3.5_LowerCarbon}\includegraphics[width=.8\textwidth,height=.25\textheight]{H03-heatmap2-v2.png}\hfill
%\end{figure}

\begin{figure}[H]
  \caption{Visual Comparison of Demand Satisfaction Rates with $H=0.3$ for $B=3M$ (left) and for $B=3.5M$ (right)}
 \begin{minipage}{0.5\linewidth}
  \centering
  \includegraphics[scale=0.6]{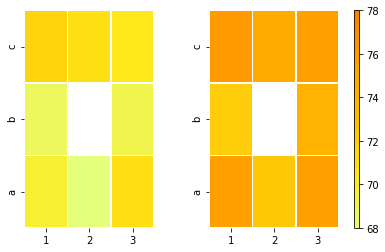}
  SROP (left) and SROP-S (right)
\label{fig:visual3_LowerCarbon}
 \end{minipage}
 \begin{minipage}{0.5\linewidth}
  \centering  \includegraphics[scale=0.6]{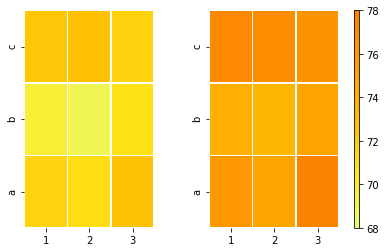}
    SROP (left) and SROP-S (right) \label{fig:visual3.5_LowerCarbon}
 \end{minipage}
\end{figure}

\iffalse
\begin{center}
\begin{tikzpicture}
\begin{axis}[
    title={Demand Rates},
    xlabel={Demand},
    ylabel={Probability Density},
    xmin=0, xmax=6,
    ymin=0, ymax=1,
    xtick={0,1,2,3,4,5},
    ytick={0,0.1,0.2,0.3,0.4,0.5,0.6,0.7,0.8,0.9},
    legend pos=north east,
    ymajorgrids=true,
    grid style=dashed,
]

\addplot[
    color=blue,
    mark=square,
    ]
    coordinates {
    (0,0.8)(1,0.15)(2,0.05)(3,0)(4,0)(5,0)
    };
    \addlegendentry{\( <4T\)}

\addplot[
    color=red,
    mark=square,
    ]
    coordinates {
    (0,0.8)(1,0.2)(2,0)(3,0)(4,0)(5,0)
    };
    \addlegendentry{\( >4T\)}
    
\end{axis}
\end{tikzpicture}
\end{center}
\fi

\end{APPENDICES}

\end{document}